\documentclass[11pt]{article}
\usepackage{bbm}

\usepackage{mathrsfs}
\usepackage{amsfonts}
\usepackage{amssymb}
\usepackage{amsfonts,amssymb, mathrsfs, amsmath,   amssymb,  theorem,  float}

\newtheorem{theorem}{Theorem}[section]

\newtheorem{cor}[theorem]{Corollary}
\newtheorem{lemma}[theorem]{Lemma}
\newtheorem{definition}[theorem]{Definition}
\newtheorem{remark}[theorem]{Remark}
\newtheorem{example}[theorem]{Example}

\floatstyle{ruled}
\newfloat{algorithm}{t}{loa}
\floatname{algorithm}{Algorithm}

\def\and{\cap}

\def\proof{{\noindent\em Proof:} }

\DeclareFontFamily{U}{fsy}{} \DeclareFontShape{U}{fsy}{m}{n}{<->s*[.
9]psyr}{} \DeclareSymbolFont{der@m}{U}{fsy}{m}{n}
\DeclareMathSymbol{\diff}{\mathord}{der@m}{182}

\newcommand{\qedd}{\hspace*{\fill}$\Box$\medskip}

\floatstyle{ruled}
\newfloat{algorithm}{t}{loa}
\floatname{algorithm}{Algorithm}

\def\X{{\mathbb{X}}}
\def\Y{{\mathbb{Y}}}
\def\U{{\mathbb{U}}}
\def\V{{\mathbb{V}}}

\def\J{{\mathcal {J}}}

\def\A{{\mathcal A}}
\def\B{{\mathcal B}}
\def\C{{\mathcal C}}

\def\P{{\mathbb P}}
\def\Q{{\mathbb Q}}

\def\bu{{\mathbf{u}}}
\def\bv{{\mathbf{v}}}

\def\prem{\hbox{\rm prem}}

\def\sat{\hbox{\rm{sat}}}
\def\asat{\hbox{\rm{asat}}}

\def\max{\hbox{\rm{max}}}

\def\deg{\hbox{\rm{deg}}}

\def\init{\hbox{\rm{I}}}
\def\ord{\hbox{\rm{ord}}}
\def\lord{\hbox{\rm{Lord}}}

\def\lead{\hbox{\rm{ld}}}

\def\dim{\hbox{\rm{dim}}}

\def\lv{\hbox{\rm{lvar}}}
\def\mod{\hbox{\rm{mod}}}

\def\ord{\hbox{\rm{ord}}}
\def\Eord{\hbox{\rm{Eord}}}

\def\den{\hbox{\rm{den}}}

\def\I{\mathcal{I}}
\def\ff{{\mathcal F}}
\def\CI{{\mathcal{I}}}

\def\Q{{\mathbb Q}}

\def\F{{\mathcal {F}}}

\def\trdeg{\hbox{\rm{tr.deg}}}

\def\dtrdeg{\hbox{$\sigma$\rm{.tr.deg}}}

\def\codim{{\rm codim}}

\def\and{\cap}

\newcounter{bean}
\def\bl{\begin{list}{Step \arabic{bean}}{\usecounter{bean}}\labelwidth=34pt}
\def\el{\end{list}}

\def\deg{{\rm deg}}
\def\lvar{{\rm lvar}}

\def\init{{\rm I}}

\def\normalization1{{\rm normalization1}}
\def\normalization{{\rm normalization}}

\def\irrfactor1{{\rm irrfactor1}}
\def\irrfactor{{\rm irrfactor}}

\def\sat{{\rm sat}}

\def\aprem{{\rm aprem}}

\begin{document}

\title{Difference Chow Form\thanks{\quad Partially
       supported by a National Key Basic Research Project of China (2011CB302400) and  by a grant from NSFC (60821002)}}
\author{Wei Li and Ying-Hong Li\\ KLMM, Academy of Mathematics and Systems Science\\
Chinese Academy of Sciences, Beijing 100190, China\\
liwei@mmrc.iss.ac.cn, liyinghong10@mails.ucas.ac.cn}
\date{}
\maketitle

\begin{abstract}
In this paper, the generic intersection theory for difference
varieties is presented. Precisely, the intersection of an
irreducible difference variety of dimension $d > 0$ and order $h$
with a generic difference hypersurface of order $s$ is shown to be
an irreducible difference variety of dimension $d-1$ and order
$h+s$. Based on the intersection theory, the difference Chow form
for an irreducible difference variety is defined. Furthermore, it is
shown that the difference Chow form of an irreducible difference
variety $V$ is transformally homogenous and has the same order as
$V$.

\end{abstract}

\section{Introduction}

Difference algebra founded by Ritt and Cohn aims to study algebraic
difference equations in a similar way that polynomial equations are
studied in algebraic geometry and differential equations are studied
in differential algebra~\cite{Cohn}. Therefore, the basic concepts
of difference algebra are similar to those of differential algebra,
which are based on those of  algebraic geometry.

The Chow form, also known as the Cayley form or the Cayley-Chow
form, is a basic concept in algebraic geometry \cite{gelfand,hodge}
and has many important applications in transcendental number theory
\cite{nes1,ph1}, elimination theory \cite{brownawell,eisenbud}, and
algebraic computational complexity \cite{Complexitychowform}.

Recently, the theory of differential Chow forms in both affine and
projective differential algebraic geometry was developed
\cite{gao-dcf,li-pdcf}. It is shown that most of the basic
properties of algebraic Chow form can be extended to its
differential counterpart \cite{gao-dcf}. Closely related with
differential Chow form, a theory of differential resultant and
sparse differential resultant was also given  \cite{gao-dcf,li-sdr}.
Furthermore, a theory of sparse difference resultants has been
developed ~\cite{li-sddr}. So it is worthwhile to generalize the
differential Chow form to its difference counterpart.

In this paper, we will study the difference Chow form for
irreducible difference varieties. We first consider the dimension
and order for the intersection of an irreducible  difference variety
by a generic difference hypersurface. Precisely, the intersection of
an irreducible difference variety of dimension $d>0$ and
order $h$ with a generic difference hypersurface of order $s$ is shown to be an
irreducible difference variety of dimension $d-1$ and order $h+s$.
Based on the intersection theory, the concept of difference Chow
form for an irreducible difference variety is defined. Furthermore,
it is shown that the difference Chow form of an irreducible
difference variety $V$ is transformally homogenous and has the same
order as $V$. The theory of characteristic set for reflexive prime difference ideals
\cite{Gao1,cohn-manifolds} plays a key role in the development of
the theory of difference Chow form.

Although both the generic intersection theorem and the basic properties of difference Chow form are similar to their differential counterparts given in \cite{gao-dcf},
some of them are quite different in terms of descriptions and proofs.
Firstly, the proof of the generic intersection theorem 
is quite different from its differential counterpart. 
In differential case, Kolchin's theory on primitive elements plays a crucial role in the proof of \cite[Theorem 3.13]{gao-dcf}.
However, the difference analogue of such theory is too weak  to be applied here.
Secondly, the definition of the difference Chow form is more subtle than the differential case 
and  the correspondence between irreducible difference varieties and the difference Chow forms may not be one-to-one as illustrated in Example~\ref{ex-chowcontrast}. 
The main reason lies in the fact that  extensions of  difference fields are much more complicated than the differential case.
Finally, the theory of difference Chow form is much more incomplete than the differential Chow form.
For instance, it lacks Poisson-type factorization formula  
and whether a theory of difference Chow variety can be developed is still open. 

The rest of the paper is organized as follows.
In section 2, we present the  basic notations and preliminary results in difference algebra.
We devote section 3 to a discussion of order and dimension for a reflexive prime difference ideal
in terms of its characteristic sets. 
Generic linear transformations as well as a generic intersection theorem on difference varieties with generic hyperplanes are then given in section 4.
And in section 5, the generic intersection theory for generic difference polynomials is established.
The difference Chow form for an irreducible difference variety is defined and its basic properties are given in section 6.
In section 7, we present the conclusion and propose several problems for further study.

\section{Preliminaries}\label{sec-preliminaries}
In this section, some notions and preliminary results in difference algebra will be given.
 For more details about difference algebra, please refer to \cite{Cohn, wibmer}.

\subsection{Difference polynomial Algebra}
A difference field $\F$ is a field with a third unitary operation
$\sigma$ satisfying that for any $a,b\in\F$,
$\sigma(a+b)=\sigma(a)+\sigma(b)$, $\sigma(ab)=\sigma(a)\sigma(b)$
and $\sigma(a)=0$ iff $a=0$.
Here, $\sigma$ is called the {\em transforming operator} of $\F$.
If $a\in\F$, $\sigma(a)$ is called the transform of $a$ and  denoted by $a^{(1)}$.
More generally, for $n\in\mathbb{Z}^{+}$, $\sigma^n(a)=\sigma^{n-1}(\sigma(a))$ is called the $n$-th transform
of $a$ and denoted by $a^{(n)}$.
And by convention, $a^{(0)}=a$.
For ease of notations, set $a^{[n]}=\{a,a^{(1)},\ldots,a^{(n)}\}$ and $a^{[\infty]}=\{a^{(i)}|i\geq0\}$.
If $\sigma^{-1}(a)$ is defined for all $a\in\F$, we say that $\F$ is inversive. A typical example of difference field is $(\Q(x),\sigma)$
with $\sigma(f(x))=f(x+1)$.

Let $S$ be a subset of a  difference field $\mathcal{G}$ which
contains $\mathcal {F}$.   We will   denote respectively by
$\mathcal {F}[S]$, $\mathcal {F}(S)$, $\mathcal {F}\{S\}$  and
$\mathcal {F}\langle S\rangle$    the smallest subring, the smallest
subfield, the smallest difference subring and the smallest
difference subfield of $\mathcal{G}$ containing both $\mathcal {F}$ and
$S$.  If we denote $\Theta(S)=\{\sigma^ka|k\geq0,a\in S\}$,  then
$\mathcal {F}\{S\}=\mathcal    {F}[\Theta(S)]$ and $\mathcal
{F}\langle S\rangle=\mathcal    {F}(\Theta(S))$.

A subset $\Sigma$ of a difference extension field $\mathcal {G}$ of
$\mathcal {F}$ is said to be {\em transformally dependent} over
$\mathcal {F}$ if the set $\Theta(\Sigma)$ is algebraically dependent over $\mathcal {F}$,
and is said to be {\em transformally independent} over $\mathcal
{F}$, or to be a family of {\em difference indeterminates} over
$\mathcal {F}$ in the contrary case.
In the case $\Sigma$ consists of only one element $\alpha$, we say that
$\alpha$ is transformally algebraic or transformally transcendental
over $\mathcal {F}$ respectively. The maximal subset $\Omega$ of
$\mathcal {G}$ which are transformally independent over $\mathcal
{F}$ is said to be a transformal transcendence basis of $\mathcal
{G}$ over $\mathcal {F}$. We use $\dtrdeg \,\mathcal {G}/\mathcal
{F}$  to denote the {\em transformal transcendence degree} of
$\mathcal {G}$ over $\mathcal {F}$, which is the cardinal number of
$\Omega$. Regarding $\mathcal {F}$ and $\mathcal {G}$ as ordinary
algebraic fields, we denote the algebraic transcendence degree of
$\mathcal {G}$ over $\mathcal {F}$ by $\trdeg\,\mathcal {G}/\mathcal
{F}$.

A homomorphism (resp. isomorphism) $\varphi$ from a difference ring $(\mathcal{R},\sigma)$ to a difference ring $(\mathcal {S},\sigma_1)$ is a {\em
difference homomorphism } (resp. {\em difference isomorphism}) if $\varphi\circ \sigma=\sigma_1\circ
\varphi$. If $\mathcal {R}_0$ is a common difference subring of
$\mathcal {R}$ and $\mathcal {S}$ and the homomorphism $\varphi$
leaves every element of $\mathcal {R}_0$ invariant, it is said to be
over $\mathcal {R}_0$.

Let $K$ be the underlying field of $\mathcal {F}$, that
is, an algebraic field consisting of the same elements as $\ff$. Let
$K(x_{1},x_{2},\ldots)$ be an overfield of
$K$ such that the $x_i$ form an algebraically
independent set over $K$, and $K^*$ be the
algebraic closure of $K(x_1,x_2,\ldots)$. Set
$\mathscr{U}$ to be the set consisting of all difference fields
 defined over subfields of $K^*$ and which are
difference overfields of $\mathcal {F}$. Then $\mathscr{U}$ is
called the universal system of difference overfields of $\mathcal
{F}$. 

Now suppose $\Y=\{y_{1},  \ldots, y_{n}\}$ is a set of
difference indeterminates over $\ff$.   The elements of $\mathcal
{F}\{\Y\}=\mathcal {F}[y_j^{(k)}:j=1,\ldots,n;k\in \mathbb{N}_0]$
are called {\em difference polynomials} over $\ff$,
 and $\mathcal {F}\{\Y\}$ itself is called the {\em difference polynomial ring } over $\ff$.
A difference polynomial ideal $\mathcal {I}$ in $\mathcal {F}\{\Y\}$
is an algebraic ideal which is closed under transforming, i.e.
$\sigma(\mathcal {I})\subset\mathcal {I}$. A difference ideal
$\mathcal {I}$ is called {\em reflexive} if $a^{(1)}\in\mathcal{I}$
implies that $a\in\mathcal{I}$.
 And a difference ideal $\mathcal {I}$ is prime if for any $a,b\in\ff\{\Y\}$, $ab\in \mathcal {I}$ implies that $a\in\mathcal {I}$ or $b\in\mathcal {I}$.
For convenience, a prime difference ideal is assumed not to be the
unit ideal in this paper.

An $n$-tuple over $\ff$ is an $n$-tuple of the form
$\textbf{a}=(a_1,\ldots,a_n)$ where the $a_i$ are selected from some
difference overfield of $\ff$. For a difference polynomial
$f\in\ff\{y_1,\ldots,y_n\}$, $\textbf{a}$ is called a difference
solution of $f$ if when substituting $y_i^{(j)}$ by $a_i^{(j)}$ in $f$,
the result is $0$, denoted by $f(\textbf{a})=0$.
For a set of difference polynomials $S\subset\ff\{\Y\}$, the {\em difference variety defined by $S$} is the set of all difference solutions
of $S$ with coordinates selected from a field of $\mathscr{U}$, denoted by $\mathbb{V}(S)$.
If we use $\textbf{a}=(a_1,\ldots,a_n)\in\mathscr{U}^n$ to mean that there exists
$\mathcal {G}\in\mathscr{U}$ such that $a_i\in\mathcal {G}$ for each
$i$, then $\mathbb{V}(S)=\{\textbf{a}\in\mathscr{U}^n|f(\textbf{a})=0, \forall \,f\in S\}$.
A difference variety  is called irreducible if it is not the union of two proper subvarieties.
For a difference variety $V$, $\mathbb{I}(V)=\{f\in\ff\{\Y\}|f(\textbf{a})=0, \forall\, \textbf{a}\in S\}$.
And $V$ is irreducible if and only if $\mathbb{I}(V)$ is a reflexive prime difference ideal.

 An $n$-tuple $\eta$ is called a {\em generic zero} of a difference ideal $\CI\subset\ff\{\Y\}$ if for any polynomial $P\in\ff\{\Y\}$, $P(\eta)=0$ iff $P\in\CI$.
It is well known that

\begin{lemma}\cite[p.77]{Cohn}\label{lm-gp}
A difference ideal possesses a generic zero if and only if it is a
reflexive prime difference ideal other than the unit ideal.
\end{lemma}


Given two $n$-tuples $\textbf{a}=(a_1,\ldots,a_n)$ and
$\bar{\textbf{a}}=(\bar{a}_1,\ldots,\bar{a}_n)$ over $\ff$,
$\bar{\textbf{a}}$ is called a specialization of $\textbf{a}$ over $\ff$,
or $\textbf{a}$ specializes to $\bar{\textbf{a}}$, if for every
difference polynomial $P\in\ff\{\Y\}$, $P(\textbf{a})=0$ implies
that $P(\bar{\textbf{a}})=0$.
A point $\eta \in \mathscr{U}^n$ is said to be {\em free from }the pure transformal transcendental
extension field $\mathcal {F}\langle U\rangle$ of $\mathcal {F}$ if $U$ are transformally independent over $\mathcal
{F}\langle \eta\rangle$. The following property about
difference specialization will be needed in this paper

\begin{lemma}\label{lm-special}
Let  $P_{i}(\U, \Y)\in \mathcal {F}\{\Y,\U\}$ $(i=1, \ldots, m)$
where $\U=(u_{1},\ldots,u_{r})$ and $\Y=(y_{1}, \ldots, y_{n})$ are
sets of difference indeterminates.
Suppose $\overline{\Y}=(\bar{y}_{1}, \ldots, \bar{y}_{n})$ is an $n$-tuple over $\ff$ that is free from $\ff\langle \U\rangle$.
If $P_{i}(\U, \overline{\Y})$ $(i=1, \ldots,
m)$ are transformally dependent over $\mathcal {F}\langle \U
\rangle$, then for any difference specialization $\U$ to
$\overline{\U}\subset\mathcal {F}$, $P_{i}(\overline{\U},\overline{\Y}) \,
(i=1, \ldots,  m)$ are transformally dependent over $\mathcal {F}$.
\end{lemma}

\proof It suffices to show the case $r=1$. Denote $u=u_1$. Since
$P_{i}(u, \overline{\Y})$ $(i=1, \ldots, m)$ are transformally dependent over
$\mathcal {F}\langle u \rangle$, there exist natural numbers $s$ and
$l$ such that $\P_i^{(k)}(u,\overline{\Y})\,(k\leq s)$ are algebraically
dependent over $\ff(u^{(k)}|k\leq s+l)$. When $u$ specializes to
$\bar{u}\in\ff$,  $u^{(k)}\,(k\geq 0)$ correspondingly
algebraically specialized to $\bar{u}^{(k)}\in\ff$. By
\cite[p.161]{wu}, $\P_i^{(k)}(\bar{u},\overline{\Y})\,(k\leq s)$ are
algebraically dependent over $\ff$. Thus, $P_{i}(\bar{u},\overline{\Y}) \,
(i=1, \ldots,  m)$ are transformally dependent over $\mathcal {F}$.
\qedd

\subsection{Characteristic set for a difference polynomial system}

A {\em ranking} $\mathscr{R}$ is a total order over $\Theta
(\Y)=\{\sigma^ky_i|1\leq i\leq n,k\geq0\}$ satisfying

1) $\sigma(\alpha)
>\alpha $ for all $\alpha \in\Theta (\Y)$ and

2) $\alpha_{1} >\alpha_{2} $ $\Longrightarrow$ $\sigma(\alpha_{1} )
>\sigma(\alpha_{2} )$ for arbitrary $\alpha_{1}, \alpha_{2}\in \Theta (\Y)$.

Below are two important kinds of rankings:

    1) {\em Elimination ranking}: \ $y_{i} > y_{j}$ $\Longrightarrow$ $\sigma^{k}  y_{i} >\sigma^{l}
    y_{j}$ for any $k, l\geq 0$.

    2) {\em Orderly ranking}: \
     $k>l$  $\Longrightarrow$  $\sigma^{k}y_{i} >\sigma^{l}
    y_{j}$,  for any $i,j \in \{1,  \ldots,n\}$.

Let $f$ be a difference polynomial in $\mathcal {F}\{\Y\}$ endowed
with a ranking $\mathscr{R}$. The {\em leader} of $f$ is the
greatest $y_j^{(k)}$ appearing effectively in $f$, denoted by
$\lead(f)$. In this case, we call $y_j$ the {\em leading variable
}of $f$, denoted by $\lvar(f)$.
 The leading coefficient of $f$  as a univariate polynomial in $\lead(f)$ is called the {\em initial} of $f$ and denoted by $\init_{f}$.
 The order of $f$ w.r.t. $y_i$, denoted by $\ord(f,y_{i})$, is defined as the greatest number $k$ such that $y_{i}^{(k)}$
appears effectively in $f$.  The {\em least order} of
$f$ w.r.t. $y_i$ is $\lord(f,y_i)=\min\{k|\deg(f,y_{i}^{(k)})>0\}$
  and the {\em effective order} of $f$ w.r.t. $y_i$ is $\Eord(f,y_i)=\ord(f,y_i)-\lord(f,y_i)$.
  And if $y_{i}$ does not appear in $f$,
then set $\ord(f,y_{i})=-\infty$ and $\Eord(f,y_i)=-\infty$.
 The {\em order} of $f$ is defined
as $\ord(f)=\max_{i}\,\ord(f,y_{i})$.

Let $f$ and $g$ be two difference polynomials in $\ff\{\Y\}$.
We say $g$ is higher than $f$, denoted by $g>f$, if 1) $\lead(g)>\lead(f)$, or 2)
$\lead(g)=\lead(f)=y_j^{(k)}$ and $\deg(g,y_j^{(k)})>\deg(f,
y_j^{(k)})$. Suppose $\lead(f)=y_j^{(k)}$. Then $g$ is said to be
{\em reduced} w.r.t. $f$ if  $\deg(g,y_j^{(k+l)})<\deg(f,y_j^{(k)})$
for each $l\geq 0$.
  A finite sequence of nonzero difference polynomials $\mathcal {A}=A_1,\ldots,A_m$  is called an
    {\em ascending chain} if  1) $m=1$ and $A_1\neq0$ or
   2) $m>1$, $A_j>A_i$ and $A_j$ is reduced
    w.r.t. $A_i$ for $1\leq i<j\leq m$.

\begin{definition}
Let $\mathcal{A}$ be an ascending chain.
Rewrite $\mathcal{A}$ in the following form
\begin{equation} \label{eq-charset}
\mathcal{A}=\left\{\begin{array}{l}A_{11},\ldots,A_{1l_1}\\ \, \cdots \\ A_{p1},\ldots,A_{pl_p} \end{array}
\right.
\end{equation} where $\lv(A_{ij})=y_{c_i}$ for $j=1,\ldots,l_i$, $c_i\ne
c_j$ for $i\ne j$, and $\ord(A_{ij},y_{c_i})<\ord(A_{ik},y_{c_i})$
for $j<k$. Then the {\em order} of $\mathcal {A}$ is defined as
$\sum_{i=1}^p \ord(A_{i1},y_{c_i})$, and the subset
$\mathbb{Y}\backslash\{y_{c_1},\ldots,y_{c_p}\}$ is called the {\em
parametric set} of $\mathcal {A}$.
\end{definition}

Set $o_{ij}=\ord(A_{ij},y_{c_i})$. For $h_1,\ldots,h_n\geq0$,
following \cite[Section 3.2]{Gao1}, we obtain the following
polynomial sequence
\[\mathcal{A}_{(h_1,\ldots,h_n)}=\left\{\begin{array}{l}A_{11},\ldots,A_{11}^{(o_{12}-o_{11}-1)},A_{12},\ldots,
A_{1l_1},\ldots,A_{1l_1}^{(\overline{h}_{c_1}-o_{1l_1})}\\
\cdots \\
A_{p1},\ldots,A_{p1}^{(o_{p2}-o_{p1}-1)},A_{p2},\ldots,
A_{pl_p},\ldots,A_{pl_p}^{(\overline{h}_{c_p}-o_{pl_p})}
\end{array} \right.\]
where $\overline{h}_{c_i}\geq \max\{h_{c_i},o_{il_i}+1\}$ is an
integer depending on $\mathcal{A}$ and the algorithm. For an
ascending chain $\mathcal {A}$ and a difference polynomial $f$, let
$\mathcal {A}_f=\mathcal {A}_{(\ord(f,y_1),\ldots,\ord(f,y_n))}$.
Then the {\em difference remainder} of $f$ w.r.t. $\mathcal {A}$ is
defined as the algebraic pseudo-remainder of $f$ w.r.t. $\mathcal
{A}_f$, that is, $\prem(f,\mathcal {A})=\aprem(f,\mathcal {A}_f)$,
which is reduced w.r.t. $\mathcal {A}$.

Let $\mathcal {A}$ be an ascending chain.
Denote $\mathbb{I}_{\mathcal {A}}$ to be the minimal multiplicative set containing the initials of elements of
$\mathcal{A}$ and their transforms. The {\em saturation ideal} of
$\A$ is defined as  $$\sat(\A)=[\mathcal
   {A}]:\mathbb{I}_{\mathcal
{A}} = \{p| \exists h\in \mathbb{I}_{\mathcal {A}}, \,{\text s. t.
}\, hp\in[A]\}.$$
The {\em algebraic saturation ideal} of $\A$ is
$\asat(\A)=(\A):\init_{\A}$, where $\init_{\A}$ is the minimal
multiplicative set containing the initials of elements of
$\mathcal{A}$.

 An ascending chain $\mathcal {C}$ contained in a difference polynomial set
 $\mathcal {S}$ is said to be a {\em characteristic set} of $\mathcal {S}$,
 if  $\mathcal {S}$ does not contain any nonzero element reduced w.r.t.
$\mathcal {C}$. If $\mathcal{C}$ is a characteristic set of a difference ideal $\CI$, then $\mathcal{C}$ reduces each element of $\CI$ to zero.
Moreover, if  $\mathcal{I}$ is a reflexive prime difference ideal, then $\mathcal {C}$ reduces to zero only the elements
of $\mathcal {I}$ and we have $\mathcal {I}=\sat(\C)$.

\begin{remark}
A set of difference polynomials $\A=\{A_1\ldots,A_m\}$ is called a
{\em difference triangular set} if the following conditions are satisfied,

1) the leaders of $A_i$ are distinct,

2) no initial of an element of $\A$ is reduced to zero by $\A$.

Similarly properties to ascending chains can be developed for difference triangular sets.
We also can define a characteristic set of a difference ideal $\mathcal{I}$ to be
a difference triangular set $\A$ contained in $\mathcal{I}$ such that $\mathcal {I}$ does not contain any nonzero element reduced w.r.t.
$\mathcal {A}$.
So we will not distinguish ascending chains and difference triangular set in this paper.
\end{remark}

\section{Dimension and order of a reflexive prime difference ideal}

Let $\CI$ be a reflexive prime difference ideal in $\mathcal
{F}\{\mathbb{Y}\}$ with a generic point $(\eta_1,\ldots,\eta_n)$.
The {\em dimension} of $\CI$ is defined as $\dtrdeg\mathcal
{F}\langle\eta _1,\ldots,\eta _n\rangle/\mathcal {F}$ and
$\codim(\mathcal {I})=n-\dim(\mathcal {I})$. A subset $\mathbb{U}$
of $\mathbb{Y}$ is called a transformal independent set modulo
$\mathcal {I}$ if $\mathcal {I} \bigcap \mathcal {F}
\{\mathbb{U}\}=\{0\}$. A maximal transformal independent set modulo
$\mathcal {I}$ is called a {\em parametric set} of $\mathcal {I}$.
Obviously, a necessary and sufficient condition for a transformal
independent set $\mathbb{U}$ becoming a parametric set of $\mathcal
{I}$ is that $|\mathbb{U}|=\dim (\mathcal {I})$. Let
$\mathbb{U}=\{y_{i_1},\ldots,y_{i_p}\}$ be a parametric set of
$\CI$, the {\em relative order of $\mathcal {I}$ w.r.t to
$\mathbb{U}$}, denoted by $\ord_\mathbb{U}(\mathcal {I})$, is
defined as  $\trdeg\,\mathcal {F}\langle\eta _1,\ldots,\eta
_n\rangle/\mathcal {F}\langle \eta_{i_1},\ldots,\eta_{i_p}\rangle$.

\begin{definition}\cite{klmp}
Let $\mathcal {I}$ be a reflexive prime difference ideal in $\mathcal {F}\{\mathbb{Y}\}$
with a generic zero $(\eta_1,\ldots,\eta_n)$.
Then there exists a numerical polynomial $\varphi_{\CI}(t)$ such that,
for all sufficiently large $t\in \mathbb{N}$, $\varphi_{\CI}(t)=\trdeg\,\mathcal {F}(\eta_1^{[t]},\ldots,\eta_n^{[t]})/\mathcal {F}$.
The polynomial $\varphi_{\CI}(t)$ is called the {\em difference dimension polynomial} of $\mathcal {I}$.
\end{definition}

\begin{lemma}\cite{klmp}\label{Dim-poly}
Let $\mathcal {I}$ be a reflexive prime difference ideal with $\dim(\mathcal {I})=d$.
Then there exists a nonnegative integer $h$ such that $\varphi_{\CI}(t)=d(t+1)+h.$
We define $h$ to be the {\em order} of $\mathcal {I}$ and denote $\ord(\mathcal {I})=h$.
Moreover, if $\mathcal {A}$ is a characteristic set of $\mathcal {I}$ w.r.t some fixed orderly ranking,
then $\ord(\mathcal {A})=\ord(\mathcal {I})$.
\end{lemma}

In \cite{cohn-manifolds}, Cohn showed that a characteristic set of $\CI$ w.r.t. an arbitrary elimination ranking
can give the information of both dimension and relative order of $\CI$.
More Precisely, if $\mathcal{A}$ is a characteristic set of $\CI$ under some elimination ranking and
$\mathbb{U}$ is the parametric set of $\mathcal{A}$,
then $\dim(\CI)=|\mathbb{U}|$ and $\ord_{\mathbb{U}}\mathcal {I}=\ord(\mathcal{A})$.
There, elimination ranking plays a crucial role in the proof.
In the following, we will show that no matter which ranking we work with,
the previous result is still true.
As an application, we also give the relation between the order and the relative order of a reflexive prime difference ideal.
Before proposing  the main result, we first give the following  lemmas.

\begin{lemma}\label{ord1}
Let $\mathcal {I}$ be a reflexive prime difference ideal in $\mathcal { F}\{\mathbb{Y}\}$ and
$\mathcal {C}$ be a characteristic set of $\mathcal {I}$.
Suppose $\mathbb{U}$ is the parametric set of $\mathcal {C}$ and set $\overline{\mathbb{Y}}=\mathbb{Y}\backslash\mathbb{U}$.
Then $\mathcal {C}$ is also a characteristic set of $\overline{\mathcal {I}}=
[\mathcal {I}]\subset\mathcal {F}\langle\mathbb{U}\rangle\{\overline{\mathbb{Y}}\}$ w.r.t. the ranking induced by the original one on $\Y$
and $\ord_\mathbb{U}(\mathcal {I})=\ord(\overline{\mathcal {I}})$.
\end{lemma}

\proof Since each nonzero $f\in\ff\{\mathbb{U}\}$ is reduced w.r.t. $\mathcal {C}$, $\CI\bigcap\ff\{\mathbb{U}\}=\{0\}$.
Hence, $\overline{\mathcal {I}}\neq[1]$.
For each $f\in \overline{\mathcal{I}}$,
there exists $h\in \mathcal {F}\{U\}$ such that $hf\in\mathcal {F}\{\mathbb{Y}\}\bigcap \overline{\mathcal {I}}=\mathcal{I}$.
If $f$ is reduced w.r.t. $\mathcal{C}$, then $hf\in\CI$ is also reduced w.r.t. $\mathcal{C}$.
 So $hf=0$ and $f=0$ follows.
 Thus, $\mathcal {C}$ is a characteristic set of $\overline{\mathcal {I}}$.
Suppose $|U|=p$ and $(U,\eta_1,\ldots,\eta_p)$ is a generic point of $\CI$.
Then $(\eta_1,\ldots,\eta_p)$ is a generic point of $\overline{\mathcal {I}}$.
By the definition of relative order, $\ord_\mathbb{U}(\mathcal {I})=\ord(\overline{\mathcal {I}})$. \qedd

\begin{lemma}\label{ord2}
 Let $\mathcal {I}$ be a reflexive prime difference ideal in $\mathcal { F}\{\mathbb{Y}\}$ and
$\mathcal {C}$ be a characteristic set of $\mathcal {I}$ w.r.t. an arbitrary ranking which has empty parametric set.
Then for sufficiently large $r\in\mathbb{N}$, the algebraic dimension of $\sat(\mathcal{C})\bigcap \mathcal {F}[\Theta_r\mathbb{Y}]$
is equal to the order of $\mathcal {C}$, where $\Theta_r\mathbb{Y}=\{y_j^{(k)}|k\leq r;j=1,\ldots,n\}$.
\end{lemma}

\proof Let $r_0=\max\{\ord(A)|A\in\mathcal {C}\}+1$ and take $r\geq
r_0$. Denote $\mathcal {C}_{r}=\mathcal {C}_{(r,\ldots,r)}$.
Firstly, since for any $f\in\sat(\mathcal{C})\bigcap \mathcal
{F}[\Theta_r\mathbb{Y}]$, $\prem(f,\mathcal{C})=\aprem(f,\mathcal
{C}_{r})=0$, it follows that $\sat(\mathcal{C})\bigcap \mathcal
{F}[\Theta_r\mathbb{Y}] =\asat(\mathcal {C}_{(r)})\bigcap \mathcal
{F}[\Theta_r\mathbb{Y}]$. By \cite[Thm 3.2]{Hubert}, the set of
non-leaders of $\mathcal {C}_{(r)}$ ia a parametric set of
$\asat(\mathcal {C}_{(r)})$. And by \cite[Lemma 3.3]{Gao1}, the set
of non-leaders of $\mathcal {C}_{(r)}$ is contained in
$\Theta_r\mathbb{Y}$ and its cardinal is equal to $\ord(\mathcal
{C})$. Thus $\dim(\asat(\mathcal {C}_{(r)})\bigcap  \mathcal
{F}[\Theta _r\mathbb{Y}])=\ord(\mathcal {C})$. Hence,
$\dim(\sat(\mathcal{C})\bigcap \mathcal
{F}[\Theta_r\mathbb{Y}])=\ord(\mathcal {C}).$ \qedd

With the above preparations, we now give the first main result in this section,
which generalizes  \cite[Theorem 4.11]{hubert} to the difference case.
\begin{theorem}\label{ord-thm}
Let $\mathcal {C}$ be a characteristic set of a reflexive prime difference ideal $\mathcal {I}\subset\mathcal {F}\{\mathbb{Y}\}$
endowed with an arbitrary ranking. Then the parametric set $\mathbb{U}$ of $\mathcal {C}$ is a parametric set of $\mathcal {I}$.
Its cardinal gives the difference dimension of $\mathcal {I}$.
Furthermore, the order of $\mathcal {I}$ relative to $\mathbb{U}$ is equal to the order of $\mathcal {C}$.
\end{theorem}

\proof Consider $\overline{\mathcal {I}}= [\mathcal
{I}]\subset\mathcal
{F}\langle\mathbb{U}\rangle\{\overline{\mathbb{Y}}\}$, where
$\overline{\mathbb{Y}}=\mathbb{Y} \backslash \mathbb{U}$. By
Lemma~\ref{ord1}, $\mathcal {C}$ is a characteristic set of
$\overline{\mathcal {I}}$ which has empty
parametric set. By Lemma~\ref{ord2},
for sufficiently large $r\in\mathbb{N}$, $\dim(\overline{\mathcal {I}}\bigcap \mathcal {F}[\Theta
_r\overline{\mathbb{Y}}])=\ord(\mathcal {C})$. By
Lemma~\ref{Dim-poly}, $\dim(\overline{\mathcal {I}}\bigcap\mathcal
{F}[\Theta_{r}\mathbb{Y}])=\dim(\overline{\mathcal
{I}})(r+1)+\ord(\overline{\mathcal {I}})$, hence
$\dim(\overline{\mathcal {I}})=0$ and $\ord(\overline{\mathcal
{I}})=\ord(\mathcal {C})$. For each $y\in
\overline{\mathbb{Y}}$, since $\overline{\mathcal {I}}\bigcap \mathcal
{F}\langle \mathbb{U}\rangle\{y\}\neq \{0\}$, $\CI\bigcap \mathcal
{F}\{\mathbb{U},y\}\neq\{0\}$. Thus $\mathbb{U}$ is a parametric
set of $\mathcal {I}$ and $\ord_\mathbb{U}\mathcal
{I}=\ord(\overline{\mathcal {I}})=\ord(\mathcal {C})$. \qedd

Apart from the trivial ideals $[0]$ and $\ff\{\Y\}$ itself,
the simplest and also most interesting ideals are reflexive prime difference ones of codimension 1.
In differential algebra, for each prime differential ideal $\CI$ of codimension $1$, there exists an irreducible differential polynomial $F$ such that $\{F\}$ is a characteristic set of $\CI$ w.r.t. any ranking.
Unlike the differential case, here even though $\I$ is of codimension one, it may happen that there are more than one difference polynomials in a characteristic set of $\I$ and
characteristic sets may be distinct for different rankings.
Nevertheless, the following lemma  shows that a uniqueness property still exists for the characteristic sets of
a reflexive prime difference ideal of codimension one under different rankings.
\begin{lemma} \cite[Lemma 2.6]{li-sddr}\label{le-char-codim1}
Let $\I$ be a reflexive prime difference ideal of codimension one in
$\ff\{\Y\}$.  The first element in any characteristic set of $\mathcal{I}$ w.r.t. any ranking, when taken irreducible,
is unique up to a factor in $\ff$.
\end{lemma}

We will end this section by proposing the following theorem,
which gives the relation between the order and the relative order of a reflexive prime difference ideal.

\begin{theorem}\label{ord-thm2}
Let $\mathcal {I}$ be a reflexive prime difference ideal in the difference
polynomial ring $\mathcal {F}\{\mathbb{Y}\}$. Then the order of
$\CI$ is equal to the maximum of all the relative orders of $\CI$,
that is, $\ord(\mathcal {I})=\max_\mathbb{U}\ord_\mathbb{U}\mathcal
{I}$, where $\mathbb{U}$ is a parametric set of $\mathcal {I}$.
\end{theorem}

\proof Let $\mathcal {C}$ be a characteristic set of $\mathcal {I}$
w.r.t some fixed orderly ranking. We claim that
$\ord_\mathbb{U}\mathcal {I}\leq \ord(\mathcal {C})$ for every
parametric set $\mathbb{U}$ of $\mathcal {I}$.

 Suppose $\mathbb{U}=\{u_1,\ldots,u_q\}$ and set $\mathbb{Y} \backslash \mathbb{U}=\{y_1,\ldots,y_p\}\,( p+q=n)$.
 Let $\mathcal {B}$ be a characteristic set of $\mathcal {I}$ w.r.t the elimination ranking $u_1\prec \ldots \prec u_q \prec y_1 \prec \ldots \prec y_p$.
Let $\eta=(\overline{u}_1,\ldots,\overline{u}_q,\overline{y}_1,\ldots,\overline{y}_p)$
be a generic zero of $\mathcal {I}$. Then for sufficiently large
 $t\in\mathbb{N}$,
\begin{eqnarray} \varphi_{\mathcal {I}}(t)
&=&\trdeg\,\mathcal {F}(\overline{u}_1^{[t]},\ldots,\overline{u}_q^{[t]},\overline{y}_1^{[t]},\ldots,\overline{y}_p^{[t]})/\mathcal {F} \nonumber \\
&=&\trdeg\,\mathcal {F}(\overline{u}_1^{[t]},\ldots,\overline{u}_q^{[t]},\overline{y}_1^{[t]},\ldots,\overline{y}_p^{[t]})
/\mathcal {F}(\overline{u}_1^{[t]},\ldots,\overline{u}_q^{[t]})+
\trdeg\mathcal {F}(\overline{u}_1^{[t]},\ldots,\overline{u}_q^{[t]})/\mathcal {F} \nonumber \\
&=&\trdeg\,\mathcal
{F}(\overline{u}_1^{[t]},\ldots,\overline{u}_q^{[t]},\overline{y}_1^{[t]},\ldots,\overline{y}_p^{[t]})/\mathcal
{F}(\overline{u}_1^{[t]},\ldots,\overline{u}_q^{[t]})+q(t+1)
\nonumber
\end{eqnarray}
Since $\varphi_{\mathcal {I}}(t)=q(t+1)+\ord(\mathcal {C})$, $\ord(\mathcal
{C})=\trdeg\mathcal
{F}(\overline{u}_1^{[t]},\ldots,\overline{u}_q^{[t]},\overline{y}_1^{[t]},\ldots,\overline{y}_p^{[t]})/\mathcal
{F}(\overline{u}_1^{[t]},\ldots,\overline{u}_q^{[t]})$. Hence,
\begin{eqnarray}
\ord_\mathbb{U}\mathcal {I}
&=&\dtrdeg\,\mathcal {F}\langle \overline{u}_1,\ldots,\overline{u}_q,\overline{y}_1,\ldots,\overline{y}_p\rangle/\mathcal {F}\langle \overline{u}_1,\ldots,\overline{u}_q\rangle \nonumber \\
&=&\trdeg\,\mathcal {F}(\overline{u}_1^{[\infty]},\ldots,\overline{u}_q^{[\infty]})
(\overline{y}_1^{[t]},\ldots,\overline{y}_p^{[t]})/\mathcal {F}(\overline{u}_1^{[\infty]},\ldots,\overline{u}_q^{[\infty]}) (t>\ord(\mathcal {B})) \nonumber \\
&\leq& \trdeg\,\mathcal {F}(\overline{u}_1^{[t]},\ldots,\overline{u}_q^{[t]},\overline{y}_1^{[t]},\ldots,\overline{y}_p^{[t]})/
\mathcal {F}(\overline{u}_1^{[t]},\ldots,\overline{u}_q^{[t]}) \nonumber \\
&=&\ord(\mathcal {C}) \nonumber
\end{eqnarray}
Set $\mathbb{U}^*$ to be the parametric set of $\mathcal {C}$, by
Lemma~\ref{Dim-poly} and Theorem~\ref{ord-thm}, $\ord(\mathcal
{I})=\ord(\mathcal {C})=\ord_{\mathbb{U}^*}\mathcal {I}$. Hence for
each parametric set $\mathbb{U},$ $\ord_\mathbb{U}(\mathcal {I})\leq
\ord(\mathcal {I})$ and there exists a parametric set $\mathbb{U}^*$ such that
$\ord_{\mathbb{U}^*}(\mathcal {I})=\ord(\mathcal {I})$. Thus, the
proof is completed. \qedd

\section{Generic linear transformations and intersection  with   generic difference hyperplanes}
In this section, we will introduce generic linear transformations and study the intersection of an irreducible
difference variety with a generic difference hyperplane. Throughout
this section, $\mathscr{U}$ stands for a fixed universal system of
difference overfields of $\mathcal {F}$ and by
$(a_1,\ldots,a_n)\in\mathscr{U}^n$, we mean that there exists
$\mathcal {G}\in\mathscr{U}$ such that $a_i\in\mathcal {G}$ for each
$i$.

\begin{definition}
Let $U=\{u_{ij}:i=1,\ldots,n;j=1,\ldots,n\}\subset\mathscr{U}$ be a
transformally independent set over $\mathcal {F}$. A generic
difference linear transformation over $\ff$ is a linear
transformation $\mathcal {T}$ from $\mathscr{U}^n$ to
$\mathscr{U}^n$ such that for every point
$\alpha=(\alpha_1,\cdots,\alpha_n)^\text{\rm T}$,
\[\mathcal {T}(\alpha)=\left(
                        \begin{array}{ccc}
                          u_{11} & \cdots & u_{1n} \\
                          \vdots & \ddots & \vdots \\
                          u_{n1} & \ldots & u_{nn} \\
                        \end{array}
                      \right)
                      \left(
                      \begin{array}{c}
                        \alpha_1 \\
                        \vdots \\
                        \alpha_n
                      \end{array}
                      \right)
\]
\end{definition}

Obviously, $\mathcal {T}$ is an inverse linear transformation and we
denote its inverse mapping to be $\mathcal {T}^{-1}$. For every
difference polynomial $P(\mathbb{Y}) \in \mathcal
{F}\{\mathbb{Y}\}$, we define $\mathcal {T}(P)=P(\mathcal
{T}^{-1}(\mathbb{Y}))\in\mathcal {F}\langle U\rangle\{\mathbb{Y}\}$.

%

\begin{lemma}\label{inter2}
If  $V$ is a difference variety over $\mathcal {F}$, then $\mathcal
{T}(V)$ is a difference variety over $\mathcal {F}\langle U\rangle$.
Furthermore, if $V$ is irreducible, then $\mathcal {T}(V)$ is also
irreducible, and $\dim(V)=\dim(\mathcal {T}(V))$,
$\ord(V)=\ord(\mathcal {T}(V))$.
\end{lemma}
\proof Suppose $V=\mathbb{V}(f_1,\ldots,f_m)$, where $f_i\in
\mathcal {F}\{\mathbb{Y}\}$. Claim: $\mathcal
{T}(V)=\mathbb{V}(\mathcal {T}(f_1),\ldots,\mathcal {T}(f_m))$. For
any $\textbf{b}\in \mathcal {T}(V)$, there exists $\textbf{a}\in V$
such that $\textbf{b}=\mathcal {T}(\textbf{a})$, then $\mathcal
{T}(f_i)(\textbf{b})=f_i(\mathcal
{T}^{-1}(\textbf{b}))=f_i(\textbf{a})=0$. Hence, $\textbf{b} \in
\mathbb{V}(\mathcal {T}(f_1),\ldots,\mathcal {T}(f_m))$. Conversely,
for any $\textbf{b}\in \mathbb{V}(\mathcal {T}(f_1),\ldots,\mathcal
{T}(f_m))$, then $\mathcal {T}(f_i)(\textbf{b})=f_i(\mathcal
{T}^{-1}(\textbf{b}))=0$. So $\mathcal {T}^{-1}(\textbf{b}) \in V$
and $\textbf{b} \in \mathcal {T}(V)$. Thus, $\mathcal
{T}(V)=\mathbb{V}(\mathcal {T}(f_1),\ldots,\mathcal {T}(f_m))$ and
$\mathcal {T}(V)$ is a difference variety over $\mathcal {F}\langle
U\rangle$.

Suppose $V$ is irreducible and $\xi$ is a generic zero of $V$ that
is free from $\mathcal {F}\langle U\rangle$. It is easy to show that
$\mathcal {T}(\xi)$ is a generic zero of $\mathcal {T}(V)$ over
$\mathcal {F}\langle U\rangle$. Indeed, for each $f\in\mathcal
{F}\langle U\rangle\{\Y\}$ satisfying $f(\mathcal {T}(\xi))=0$,
$\mathcal {T}^{-1}(f)(\xi)=0$. Since $\xi$ is  free from $\mathcal
{F}\langle U\rangle$,  for each $\textbf{a}\in V$, $\mathcal
{T}^{-1}(f)(\textbf{a})=0=f(\mathcal {T}(\textbf{a}))$. Thus,
$f|_{\mathcal {T}(V)}\equiv0$ and it follows that $\mathcal {T}(V)$
is irreducible.

Suppose $\dim(V)=d$, $\ord(V)=h$,  then for  sufficiently large
$t\in\mathbb{N}$, $ \varphi_{\mathcal {T}(V)/\mathcal {F}\langle
U\rangle}(t) =\trdeg\,\mathcal {F}\langle U\rangle(\mathcal
{T}(\xi)^{[t]})/\mathcal {F}\langle U\rangle =\trdeg\,\mathcal
{F}\langle U\rangle(\xi^{[t]})/\mathcal {F}\langle U\rangle
=\trdeg\,\mathcal {F}(\xi^{[t]})/\mathcal {F} =d(t+1)+h. $ Hence,
$\dim(V)=\dim(\mathcal {T}(V))$, $\ord(V)=\ord(\mathcal {T}(V))$.
\qedd

\begin{definition}
A generic difference hyperplane is the difference variety defined
by  $u_0+u_1y_1+\cdots+u_ny_n=0$, where the $u_i\in\mathscr{U}$ are
transformally independent over $\ff$.
\end{definition}

The following theorem gives the main result of this section,
which generalize an interesting theorem \cite[ p. 54, Theorem I]{hodge} in algebraic geometry to he difference case.
\begin{theorem}\label{inter-thm}
Let V be an irreducible difference variety over $\mathcal {F}$ with
dimension $d>0$ and order h. Let $\mathcal
{L}:u_1y_1+\cdots+u_ny_n-u_0=0$ be a generic difference hyperplane.
Then $V\bigcap L$ is an irreducible difference variety over
$\mathcal {F}\langle u_0,u_1,\ldots,u_n\rangle$ with dimension $d-1$
and order h.
\end{theorem}

\proof Consider the following generic difference linear
transformation $\mathcal {T}: \mathscr{U}^n\rightarrow
\mathscr{U}^n$ over $\ff$: $\forall\,
\alpha=(\alpha_1,\cdots,\alpha_n)^{\text{\rm T}}$,
$$\mathcal {T}(\alpha)=\left(
                        \begin{array}{ccc}
                          u_1 & \ldots & u_n \\
                          v_{1} & \ldots & v_{n}\\
                          \vdots & \ddots & \vdots \\
                          w_{1} & \ldots & w_{n} \\
                        \end{array}
                      \right)
                      \left(
                      \begin{array}{c}
                        \alpha_1 \\
                        \alpha_2\\
                        \vdots \\
                        \alpha_n
                      \end{array}
                      \right),
$$
where $U=\{u_i,v_j,\ldots,w_k\}$ is a transformally independent set
over $\mathcal {F}$.

Then $\mathcal {T}(\mathcal {L})= y_1-u_0=0$.  By lemma~\ref{inter2}
$\mathcal {T}(V)$ is an irreducible difference variety over
$\mathcal {F}\langle U\rangle$ with $\dim(\mathcal {T}(V))=d$ and
$\ord(\mathcal {T}(V))=h$. Suppose
\[\mathcal{A}=\left\{\begin{array}{l}A_{11},\ldots,A_{1l_1}\\ \cdots \\ A_{n-d,1},\ldots,A_{n-d,l_{n-d}} \end{array}
\right.\] is a difference characteristic set of $\mathcal {T}(V)$
w.r.t. some orderly ranking $\mathscr{R}$, where
$\lv(A_{ij})=y_{c_i}$ for $j=1,\ldots,l_i$ and
$\ord(A_{ij},y_{c_i})<\ord(A_{il},y_{c_i})$ for all $j<l$. By
interchanging the rows of the  matrix of $\mathcal {T}$ when
necessary, suppose $y_1$ lies in  the parametric set of
$\mathcal{A}$.

In each $A_{ij}$, replace $y_1$ by $u_0$ and denote it by $B_{ij}$.
Set $B_0=y_1-u_0$ and
\[\mathcal{B}=\left\{\begin{array}{l} B_0\\ B_{1,1},\ldots,B_{1,l_1}\\ \cdots \\B_{n-d,1},\ldots, B_{n-d,l_{n-d}} \end{array}
\right.\] We claim that $[\sat(\mathcal {A}),B_0]$ is a reflexive
prime difference ideal over $\ff\langle u_0,U\rangle$ and $\mathcal
{B}$ is a characteristic set of it w.r.t. $\mathscr{R}$. Then
$\mathcal {T}(V\cap L)=\mathcal {T}(V)\bigcap \mathcal {T}(L)$ is an
irreducible difference variety over $\mathcal {F}\langle
u_0,U\rangle$ with dimension $d-1$ and order $h$. Since $\mathcal
{T}$ is a inverse linear difference transformation, $V\cap\mathcal
{L}$ is an irreducible difference variety over $\mathcal {F}\langle
u_0,u_1,\ldots,u_n\rangle$ with dimension $d-1$ and order $h$. Thus
it suffices to prove the above claim.

Suppose $\zeta=(u_0,y_2,\ldots,y_d,\eta_{d+1},\ldots,\eta_n)$ is a
generic point of $\sat(\mathcal {A})\subset\mathcal {F}\langle
U\rangle\{\mathbb{Y}\}$. Let $\mathbb{I}(\zeta)$ be the difference
polynomial ideal in $\ff\langle U,u_0\rangle\{\mathbb{Y}\}$ with
$\zeta$ as a generic point. Clearly, $[\sat(\mathcal
{A}),B_0]\subset\mathbb{I}(\zeta)$. Conversely, for any $f \in
\mathbb{I}(\zeta)$, there exists $M(u_0)\in\ff\{u_0\}$ such that
$M(u_0)f\in \mathcal {F}\langle U\rangle\{u_0,\mathbb{Y}\}$. Then
$M(u_0)f(\zeta)=0$. Let $f_1=\prem(M(u_0)f,y_1-u_0)$, i.e.
$M(u_0)f\equiv f_1, \,\mod\, [y_1-u_0]$, then $f_1\in\mathcal
{F}\langle U\rangle\{u_0,\mathbb{Y}\}$ is free from $y_1$. On the
one hand, replace $u_0$ by $y_1$ in $f_1$ and denote the obtained
polynomial by $\widetilde{f_1}$, then $\widetilde{f_1}\in \mathcal
{F}\langle U\rangle \{\mathbb{Y}\}$ vanishes at $\zeta$ and
$\widetilde{f_1}-f_1 \in [y_1-u_0]$. Hence $\widetilde{f_1}\in
\sat(\mathcal {A})$ and $f_1\in [\sat(\mathcal {A}),B_0]$. Thus, $f
\in [\sat(\mathcal {A}),B_0]$. So $[\sat(\mathcal
{A}),B_0]=\mathbb{I}(\zeta)$  is a reflexive prime difference ideal.
On the other hand, let $r=\prem(f_1,\mathcal {B})$ and
$\widetilde{r}$ be obtained by replacing $u_0$ with $y_1$ in $r$.
Then $\widetilde{r}\in\sat(\mathcal {A})$ is reduced w.r.t.
$\mathcal {A}$. Thus, $\widetilde{r}=0$ and $r=0$ follows. That is,
$\mathcal {B}$ reduces all element in $\mathbb{I}(\zeta)$ to zero.
Since $\mathcal {B}\subset \mathbb{I}(\zeta) =[\sat(\mathcal
{A}),B_0]$, $\mathcal {B}$ is a characteristic set of $[\sat(\mathcal
{A}),B_0]$ w.r.t. $\mathscr{R}$. \qedd

\section{Intersection theory for generic difference polynomials }
In this section,  we will develop an intersection theory for generic difference polynomials,
which generalizes \cite[Theorem 1.1]{gao-dcf} in differential algebra to the difference case.

\begin{definition}
\label{def-genericpol}
Let $\mathbbm{m}_{s,r}$ be the set of all difference monomials in
$\ff\{\Y\}$ of  order $\leq s$ and degree $\leq r\,(r>0)$. Let
$\mathbb{U}=\{u_m\}_{m\in \mathbbm{m}_{s,r}}$ be a set of  difference indeterminates over $\ff$. Then,
$$f=\sum_{m\in \mathbbm{m}_{s,r}}u_m m$$ is called a {\em generic difference polynomial} of order $s$ and degree $r$.
A {\em generic difference hypersurface} is the set of zeros of a
generic difference polynomial.
\end{definition}

Our first goal is to show
that by adding a generic difference polynomial to a reflexive prime difference ideal, the
new ideal is still reflexive prime and its dimension will decrease by one.

\begin{theorem} \label{th-intersectionsurface-dim}
Let $\CI$ be a reflexive prime difference ideal in $\ff\{\Y\}$ of dimension $d$.
Let $\P$ be a generic difference polynomial with coefficients $\bu$.
If $d>0$, then $[\CI,\P]\subset\ff\langle\bu\rangle\{\Y\}$  is a reflexive prime difference ideal of dimension $d-1$.
And if $d=0$, then  $[\CI,\P]=\ff\langle\bu\rangle\{\Y\}$.
\end{theorem}
\proof Suppose $u_0$ is the degree zero term of $\P$.
Denote $\widetilde{\P}=\P-u_0$ and $\tilde{\bu}=\bu\backslash\{u_0\}$.
Let $\xi$ be a generic point of $\CI$ over $\ff$ that is free from $\bu.$
Let $\mathcal {J}=[\CI,\P]\subset\ff\langle\tilde{\bu}\rangle\{\Y,u_0\}$.
It is easy to show that $(\xi,-\widetilde{\P}(\xi))$ is a generic point of $\mathcal {J}$.
Thus, $\mathcal {J}$ is a reflexive prime difference ideal with
 $\dim(\mathcal {J})=\dtrdeg\,\ff\langle\tilde{\bu}\rangle\langle\xi,-\widetilde{\P}(\xi)\rangle/\ff\langle\tilde{\bu}\rangle
 =\dtrdeg\,\ff\langle\tilde{\bu}\rangle\langle\xi\rangle/\ff\langle\tilde{\bu}\rangle=d$.
If $d=0$, then  $\mathcal {J}\cap\ff\langle\tilde{\bu}\rangle\{u_0\}\neq\{0\}$.
Thus, in $\ff\langle\bu\rangle\{\Y\}$,  $[\CI,\P]=\ff\langle\bu\rangle\{\Y\}$.

It remains to show the case $d>0$.
 Without loss of generality, suppose $\{y_1,\ldots,y_d\}$ is a parametric set of $\CI$.
We claim that $\{y_1,\ldots,y_{d-1},u_0\}$ is a parametric set of $\mathcal {J}$ over $\ff\langle\tilde{\bu}\rangle$.
Suppose the contrary. Then $\xi_1,\ldots,\xi_{d-1},-\widetilde{\P}(\xi)$ are transformally dependent over $\ff\langle\tilde{\bu}\rangle$.
Now specialize the coefficient of $y_k$ in $\P$ to $-1$ and all the other $u\in\tilde{\bu}$ to zero,
then by Lemma~\ref{lm-special}, $\xi_1,\ldots,\xi_d$ are transformally dependent over $\ff$, a contradiction.
So $\mathcal {J}\cap\ff\langle\tilde{\bu}\rangle\{ y_1,\ldots,y_{d-1},u_0\}=\{0\}$.
Thus, $[\mathcal {J}]\cdot\ff\langle\bu\rangle\{\Y\}\neq[1]$ is a reflexive prime difference ideal and $[\mathcal {I},\P]\cap\ff\langle\bu\rangle\{ y_1,\ldots,y_{d-1}\}=\{0\}$.
For each $y_k\,(k\geq d)$, since $\mathcal {J}\cap\ff\langle\tilde{\bu}\rangle\{ y_1,\ldots,y_{d-1},y_k,u_0\}\neq\{0\}$,
$[\mathcal {I},\P]\cap\ff\langle\bu\rangle\{ y_1,\ldots,y_{d-1},y_k\}\neq\{0\}.$
 Hence, $[\CI,\P]\subset\ff\langle\bu\rangle\{\Y\}$  is a reflexive prime difference ideal of dimension $d-1$.
\qedd

Next, we consider the order of the intersection of an irreducible difference variety by a generic difference hypersurface.
Before proving the main result, we need two lemmas.

\begin{lemma} \label{le-elim-charset}
Let $\CI$ be a reflexive prime difference ideal in $\ff\{u_1,\ldots,u_q,y_1,\ldots,y_p\}$ and
 \[\mathcal{A}=\left\{\begin{array}{l}A_{11},\ldots,A_{1l_1}\\A_{21},\ldots,A_{2l_2}\\ \,\cdots \\ A_{p1},\ldots,A_{pl_p} \end{array}
\right.\] be a characteristic set of $\CI$ w.r.t. the elimination ranking $u_1\prec\cdots\prec u_q\prec y_1\prec \cdots\prec y_p$ with $\lv(A_{ij})=y_i$.
Suppose $f\in\CI$  is reduced w.r.t. $A_{21},\ldots,A_{pl_p}$.
Rewrite $f$ in the following form:  $f=\sum_{\phi}f_{\phi}(u_1,\ldots,u_q,y_1)\phi(y_2,\ldots,y_p)$,
where $\phi$ rangs over all distinct difference monomials appearing effectively  in $f$ and $f_{\phi}\in\ff\{u_1,\ldots,u_q,y_1\}$.
Then for each $\phi$, $f_{\phi}(u_1,\ldots,u_q,y_1)\in\CI$.
\end{lemma}
\proof
Denote $\mathcal{B}=A_{11},\ldots,A_{1l_1}$ and $\mathcal{B}_{f}=B_{1},\ldots,B_s$.
Since $f\in\CI$ is reduced w.r.t. $A_{21},\ldots,A_{pl_p}$, the difference remainder of $f$ w.r.t. $\mathcal{B}$ is zero.
Thus, $\aprem(f,\mathcal {B}_f)=0$.
Suppose $\lead(B_i)=y_1^{(o_1+i-1)}\,(i=1,\ldots,s)$.
Now we proceed to construct an algebraic triangular set $\mathcal{C}=C_1,\ldots,C_s$ contained in $\CI$ such that 1) $\lead(C_i)=\lead(B_i)$, 2) $\init_{C_i}\in\ff\{u_1,\ldots,u_q\}[y_1^{[o_1-1]}]$ and 3) $\aprem(f,\mathcal {C})=0$.
Set $C_1=B_1$. For $i=2$, if $\ord(\init_{B_2},y_1)=o_1-1$, then set $C_2=B_2$.
Otherwise, $\ord(\init_{B_2},y_1)=o_1$.
Let $R$ be  the Sylvester resultant of $\init_{B_2}$ and $B_1$ w.r.t. $y_1^{(o_1)}$.
Since $\mathcal{B}_{f}$ is a regular chain\cite[Theorem 4.1]{Gao1},  $R\neq0$  and there exist polynomials $D_1,D_2$ such that $R=D_1B_1+D_2\init_{B_2}$.
Let $C_2=\aprem(D_2B_2,B_1)$.
Clearly, $C_2\in\CI$, $\lead(C_2)=y_1^{(o_1+1)}$ and $\init_{C_2}=R\in\ff\{u_1,\ldots,u_q\}[y_1^{[o_1-1]}]$.
Similarly in this way, $\mathcal{C}=C_1,\ldots,C_s$ can be constructed.

For each $\phi$, let $r_\phi=\aprem(f_\phi,\mathcal{C})$.
 Then there exist integers $l_{\phi i}$ s.t. $\prod_{i=1}^s(\init(C_i))^{l_{\phi i}}f_\phi\equiv r_\phi\,\mod\,(\mathcal{C})$.
Let $l=\max_{\phi,i}\{l_{\phi i}\}$.
 Then  $\sum_{\phi}\prod_{i=1}^s(\init(C_i))^{l-l_{\phi i}}r_{\phi}\phi(y_2,\ldots,y_p)$ belongs to $\CI$ and is reduced w.r.t. $\mathcal{B}$.
 Thus, for each $\phi$, $r_\phi=0$ and $f_\phi\in\CI$ follows.  \qedd

\begin{lemma} \label{le-order+1}
Let $\mathcal {S}$ be a system of difference polynomials in $\mathcal {F}\{\Y\}$.
Suppose $\V(S)$ has an irreducible component $V$ of dimension $d$ and order $h$ with $\mathbb{I}(V)\cap\ff\{y_1\}=\{0\}$.
Let $\mathcal{\overline{S}}$ be obtained from $\mathcal {S}$ by replacing $y_{1}^{(k)}$ by
$y_{1}^{(k+1)}(k=0, 1,\ldots)$ in all of the polynomials in $\mathcal {S}$.
Then the variety of $\overline{\mathcal {S}}$ has a component $\overline{V}$ of dimension $d$ and order $h_1$ such that
$h\leq h_1\leq h+1$.
\end{lemma}
\proof Let $\eta=(\eta_1,\ldots,\eta_n)$ be a generic point of $V$.
Since $\I(V)\cap\ff\{y_1\}=\{0\}$, $\eta_1$ is transformally transcendental over $\ff$.
So $\mathcal{I}=[z^{(1)}-\eta_1]\subset\ff\langle\eta\rangle\{z\}$ is a reflexive prime difference ideal of dimension $0$ and order $1$.
Let $\zeta$ be a generic point of $\CI$.
Clearly, $(\zeta,\eta_2,\ldots,\eta_n)$  is a difference solution of $\mathcal{\overline{S}}$.

Suppose $(\zeta,\eta_2,\ldots,\eta_n)$ lies in a component $\overline{V}$ of $\overline{\mathcal {S}}$, which has a
generic point $\xi=(\xi_{1},\ldots,  \xi_{n})$. Then
$(\xi_{1}, \xi_2, \ldots,  \xi_{n})$ specializes to $(\zeta,\eta_{2}, \ldots, \eta_{n})$, and
 $(\xi_{1}^{(1)}, \xi_2,  \ldots,\xi_{n})$ specializes to $(\eta_{1},\eta_{2},
\ldots, \eta_{n})$ correspondingly. Since $(\eta_{1},\eta_{2},
\ldots, \eta_{n})$ is a generic point of $V$ and  $(\xi_{1}^{(1)}, \xi_{2}, \ldots,\xi_{n})$ is a zero of
$\mathcal {S}$,  $(\xi_{1}^{(1)}, \xi_{2}, \ldots,\xi_{n})$ is a generic point of $V$.
So for sufficiently large $t\in\mathbb{N}$, $\trdeg\,\ff(\xi_1^{(1)},\ldots,\xi_1^{(t+1)},\xi_2^{[t]},\ldots,\xi_n^{[t]})/\ff=d(t+1)+h$.
Since $\trdeg\,\ff\langle\eta\rangle(\zeta)/\ff\langle\eta\rangle=1$,
$\trdeg\,\ff\langle\xi_1^{(1)},\xi_2,\ldots,\xi_n\rangle(\xi_1)/\ff\langle\xi_1^{(1)},\xi_2,\ldots,\xi_n\rangle=1$.
Thus,
\begin{eqnarray}
\varphi_{\overline{V}}(t)&=&\trdeg\,\mathcal {F}(\xi_1^{[t]},\ldots,\xi_n^{[t]})/\mathcal {F} \nonumber\\
&=& \trdeg\,\mathcal {F}(\xi_1^{[t+1]},\ldots,\xi_n^{[t]})/\mathcal {F}-\trdeg\,\mathcal {F}(\xi^{[t]})(\xi^{(t+1)})/\mathcal {F}(\xi^{[t]}) \nonumber\\
&=& d(t+1)+h+1-\trdeg\,\mathcal {F}(\xi^{[t]})(\xi^{(t+1)})/\mathcal {F}(\xi^{[t]}). \nonumber
\end{eqnarray}
Consequently, $\varphi_{\overline{V}}(t)=d(t+1)+h_1$ where $h\leq h_1\leq h+1$.
By Lemma~\ref{Dim-poly}, the proof is completed.
   \qedd

With the above preparations, we now propose the main theorem in this section.
\begin{theorem} \label{th-intersectionord}
Let $\CI$ be a reflexive prime difference ideal in $\ff\{\Y\}$ of dimension $d>0$ and order $h$.
Let $\P$ be a generic difference polynomial with coefficient set $\bu$.
Then $[\CI,\P]\subset\ff\langle\bu\rangle\{\Y\}$  is a reflexive prime difference ideal of dimension $d-1$ and order $h+s$
\end{theorem}
\proof Let $\CI_1=[\CI,\P]\subset\ff\langle\bu\rangle\{\Y\}$.
By Theorem~\ref{th-intersectionsurface-dim},  $\mathcal {I}_1$ is a reflexive prime difference ideal of dimension $d-1$.
We only need to show that the order of $\mathcal {I}_1$ is $h+s$.

Let $\mathscr{A}$ be a characteristic set of $\mathcal {I}$ w.r.t.
some orderly ranking $\mathscr{R}$ with $y_{1}, \ldots, y_{d}$ as a
parametric set.
By Theorem \ref{Dim-poly}, $\ord(\mathscr{A})=h$.
Let $u_0$ be the degree zero term of $\P$ and $\widetilde{\bu}=\bu\backslash\{u_0\}$.
Let $\mathcal{J}=[\mathcal {I},\P]\subset\mathcal {F}\langle
\widetilde{\bu}\rangle\{\Y, u_{0}\}$.
By the proof of Theorem~\ref{th-intersectionsurface-dim}, $\mathcal {J}$ is a reflexive prime difference ideal of
dimension $d$. Clearly, $\mathcal {I}_1\bigcap\mathcal {F}\langle
\widetilde{\bu}\rangle\{\Y, u_{0}\}=\mathcal{J}$. So any characteristic set of $\mathcal {I}_1$, by
clearing denominators in $\ff\langle
\widetilde{\bu}\rangle\{u_0\}$ when necessary, is a characteristic
set of $\mathcal {J}$ with $u_{0}$ in the parametric set.
By Theorem \ref{ord-thm2},  we have
$\ord(\mathcal {I}_1)\leq \ord(\mathcal {J})$.

We claim that $\ord(\mathcal {J})\le h+s$.
As a consequence,  $\ord(\mathcal {I}_1)\leq h+s$.
To prove this claim, let $\mathcal {J}^{(i)}=[\mathcal {I},
u_{0}^{(i)}+\widetilde{\P}]\,(i=0,\ldots,s)$ in  $\mathcal {F}\langle
\widetilde{\bu} \rangle\{\Y, u_{0}\}.$
Similarly to the proof of Theorem~\ref{th-intersectionsurface-dim},
we can show that $\mathcal{J}^{(i)}$ is a reflexive prime difference ideal of dimension $d$. Let
$F$ be the difference remainder of $u_{0}^{(s)}+\widetilde{\P}$ w.r.t.
$\mathcal {A}$ under the ranking $\mathscr{R}$. Clearly,
$\ord(F,u_{0})=s.$ It is obvious that for some orderly
ranking, $\{\mathcal{A}, F \}$ is a characteristic set of
$\mathcal {J}^{(s)}$ with $y_1,\ldots,y_d$ as a parametric set.
So $\ord(\mathcal {J}^{(s)})=h+s$.
%
%
Using Lemma~\ref{le-order+1} $s$ times, we have $\ord(\mathcal {J})\le\ord(\mathcal{J}^{(1)})\le\cdots\le \ord(\mathcal
{J}^{(s)})=h+s$.

Now, it remains to prove $\ord(\mathcal {I}_1)\geq h+s$.
Let $\P=u_0+\sum_{i=1}^n \sum_{j=0}^s u_{ij} y_{i}^{(j)}+T$, where $T$ is the nonlinear part of $\P$.
Let $w=u_{0}+\sum_{i=1}^d \sum_{j=0}^s u_{ij} y_{i}^{(j)}$ be a new difference indeterminate.
Let $\bu_G$ be the set of coefficients of $G=w+\sum_{i=d+1}^n \sum_{j=0}^s u_{ij} y_{i}^{(j)}+T$ regarded as
a difference polynomial in $w$ and $\Y$.
We denote $\mathcal{F}_1=\mathcal {F}\langle\bu_G\rangle$.
It is easy to show that $\mathcal{J}_1=[\mathcal {I},G] \subset\mathcal {F}_1\{y_{1}, \ldots, y_{n},w\}$ is a reflexive prime difference ideal with a generic point $(\xi,-\sum_{i=d+1}^n \sum_{j=0}^s u_{ij} \xi_{i}^{(j)}-T(\xi))$,
where $\xi=(\xi_1,\ldots,\xi_n)$ is a generic point of $\mathcal{I}$.
So $y_1,\ldots,y_d$ is a parametric set of $\mathcal{J}_1$ and $\ord_{y_1,\ldots,y_d}\mathcal {J}_{1}=\ord_{y_1,\ldots,y_d}\mathcal {I}=h$.
Let \[\mathcal{B}=\left\{\begin{array}{l} R(y_1,\ldots,y_d,w),R_1(y_1,\ldots,y_d,w),\ldots,R_l(y_1,\ldots,y_d,w)\\ B_{11}(y_1,\ldots,y_d,w,y_{d+1}),\ldots,B_{1l_1}(y_1,\ldots,y_d,w,y_{d+1})
\\ \quad\cdots\cdots
\\B_{n-d,1}(y_1,\ldots,y_d,w,\ldots,y_n),\ldots, B_{n-d,l_{n-d}}(y_1,\ldots,y_d,w,\ldots,y_n) \end{array}
\right.\]
be a characteristic set of $\mathcal {J}_{1}$ w.r.t. the elimination ranking $y_1\prec \cdots\prec y_d\prec w\prec y_{d+1} \prec\cdots\prec y_n$.
Then by Theorem~\ref{ord-thm}, $\ord(\mathcal{B})=\ord_{y_1,\ldots,y_d}\mathcal {J}_{1}=h$.

Let $\bu_d = \{u_{ij}:\, i=1, \ldots, d; j=0,\ldots,s\}$.
Then $[\mathcal{J}_1]\subset\mathcal {F}_1\langle \bu_d \rangle\{w, y_{1}, \ldots,    y_{n}\}$ is also a reflexive prime difference ideal with $\mathcal{B}$ as a characteristic set w.r.t. the
elimination ranking $y_{1} \prec \ldots \prec y_{d}\prec w\prec
y_{d+1}\prec \ldots\prec y_{n}$. Let

\[ \begin{array}{ccc} \phi:\mathcal {F}_1\langle \bu_d
\rangle\{y_{1}, \ldots, y_{n},w\} &\longrightarrow &\mathcal
{F}_1\langle \bu_d \rangle\{y_{1},\ldots, y_{n}, u_{0}\}
\\ w & \quad & u_{0}+\sum_{i=1}^d \sum_{j=0}^s u_{ij}
y_{i}^{(j)}\\ y_{i}& \quad & y_{i} \end{array} \]

\noindent be a difference homomorphism over $\mathcal {F}_1\langle
\bu_d \rangle$. Clearly,  this is a difference isomorphism which
maps $[\mathcal {J}_{1}]$ to $\mathcal {J}$. It is obvious that
$\phi(R), \phi( R_1),\ldots,\phi(R_l),\phi(B_{11}), \ldots,
\phi(B_{n-d,l_{n-d}})$ is a characteristic set of $\mathcal{J}$ w.r.t. the elimination
ranking $y_{1}\prec \ldots \prec y_{d} \prec u_{0}\prec y_{d+1}\prec
\ldots\prec y_{n}$ and $\lead(\phi(B_{ij}))=\lead(B_{ij})\,(i=1, \ldots,
n-d;j=1,\ldots,l_{i})$.
We claim that $\ord(\phi(R), y_{1})\geq \ord(R,w)+s$.
Denote $\ord(R,w)=o$.
If $\ord(R,y_{1})\geq o+s$, rewrite $R$ in the form $R=\sum_{\psi_{\nu}(w)\neq
1} p_{\nu}(y_{1}, \ldots, y_{d})\psi_{\nu}(w)+p(y_{1}, \ldots,
y_{d})$ where $\psi_{\nu}(w)$ are monomials in $w$ and its
transforms.
Then
\begin{eqnarray}
\phi(R)&=& \sum_{\psi_{\nu}\neq 1} p_{\nu}(y_{1},    \ldots,
y_{d})\psi_{\nu}(u_{0}+\sum_{i=1}^d
\sum_{j=0}^s u_{ij}y_{i}^{(j)})+p(y_{1},    \ldots,    y_{d}) \nonumber \\
&=&\sum_{\psi_{\nu}\neq
1}p_{\nu}(y_{1},    \ldots,    y_{d})\psi_{\nu}(u_{0})+p(y_{1},    \ldots,    y_{d}) \nonumber \\
& &+\hbox{terms involving}\, u_{ij}(i=1,\ldots,
d;j=0,\ldots,s) \, \hbox{ and their transforms.} \nonumber
\end{eqnarray}
\noindent Clearly, in this case we have $\ord(\phi(R), y_{1})\geq
\max\{\ord(p_{\nu}, y_{1}),\ord(p, y_{1})\}=\ord(R,$ $ y_{1})\geq o+s$.
If $\ord(R ,y_{1})<o+s$, rewrite $R$ as a polynomial in $w^{(o)}$,
that is, $R=I_{l}(w^{(o)})^l+I_{l-1}(w^{(o)})^{l-1}+\cdots+I_{0}$.
Then $\phi(R)=\phi(I_{l})[(u_{0}^{(o)}+\sum_{i=1}^d \sum_{j=0}^s
u_{ij}^{(o)}y_{i}^{(o+j)})]^l+\phi(I_{l-1})[(u_{0}^{(o)}+\sum_{i=1}^d
\sum_{j=0}^s u_{ij}^{(o)}y_{i}^{(o+j)})]^{l-1}+\cdots+\phi(I_{0})$. Since
$\ord(\phi(I_k),y_1)<o+s$ ($k=0,\ldots,l$), we have exactly
$\ord(\phi(R), y_{1})=o+s$.  Thus, consider the two cases together,
$\ord(\phi(R), y_{1})\geq \ord(R,w)+s$.

Since $\mathcal {J} \cap \mathcal
{F}_1\langle\bu_d\rangle\{y_{1}, \ldots, y_{d},u_0\}$ is reflexive  prime difference ideal of codimension $1$,
by Lemma \ref{le-char-codim1},
$\phi(R)$ can serve as the first difference polynomial in a characteristic set of $\mathcal {J} \cap \mathcal
{F}_1\langle\bu_d\rangle\{y_{1}, \ldots, y_{d},u_0\}$ w.r.t. any ranking.
Suppose $\phi(R),\widetilde{R}_1,\ldots,\widetilde{R}_{\widetilde{l}}$ is a characteristic set of $\mathcal {J} \cap \mathcal
{F}_1\langle\bu_d\rangle\{y_{1}, \ldots, y_{d},u_0\}$ w.r.t. the elimination ranking $u_{0}\prec y_{2}\prec
\ldots \prec y_{d} \prec y_{1}$.
By lemma~\ref{le-elim-charset}, $\phi(R),\widetilde{R}_1,\ldots,\widetilde{R}_{\widetilde{l}},\phi(B_{11}),\ldots,\phi(B_{n-d,l_{n-d}})$
is a characteristic  set of $\mathcal{J}$ w.r.t. the elimination ranking
$u_{0}\prec y_{2}\prec \ldots \prec y_{d} \prec y_{1}\prec y_{d+1}
\prec \ldots \prec y_{n}$,     thus
a characteristic set of $\mathcal {I}_1$.
By Theorem~\ref{ord-thm2}, $\ord(\mathcal{I}_1)\geq\ord_{y_{2}, \ldots,y_{d}}\mathcal{I}_1=\ord(\phi(R),y_1)+\sum_{i=1}^{n-d}\ord(\phi(B_{i1}),y_{d+i})\geq \ord(R,w)+s+\sum_{i=1}^{n-d}\ord(B_{i1},y_{d+i})=\ord(\mathcal{B})+s=h+s$.
Thus, the order of $\mathcal{I}_1$ is $h+s$. \qedd

As a corollary, we give the dimension theorem for generic
difference polynomials.

\begin{theorem}
Let $f_{1},  \ldots, f_{r}\, (r\le n)$ be independent generic
difference polynomials with each $ f_{i}$ of order $s_{i}$. Then
$[f_{1}, \ldots, f_{r}]$ is a reflexive prime difference ideal of dimension $n-r$ and
order $\sum_{i=1}^r s_{i}$ over
$\ff\langle\bu_{f_1},\ldots,\bu_{f_r}\rangle$.
\end{theorem}
\proof We will prove the theorem  by induction on $r$. Let
$\CI=[0]\subset\ff\{\Y\}$. Clearly, $I$ is a reflexive prime difference ideal of
dimension $n$ and order $0$.
For $r=1$, by Theorem~\ref{th-intersectionord}, $[f_1]=[\CI,f_1]$ is a reflexive prime difference ideal of dimension $n-1$ and order $s_1$.
So the assertion holds for $r=1$.
Now suppose the assertion holds for $r-1$, we now prove it for $r$.
By the hypothesis, $\CI_{r-1}=[f_1,\ldots,f_{r-1}]$ is a reflexive prime difference ideal of dimension $n-r+1$ and order $\sum_{i=1}^{r-1}
s_{i}$ over $\ff\langle \bu_{f_1},\ldots,\bu_{f_{r-1}}\rangle$.
Since $f_1,\ldots,f_r$ are independent generic difference polynomials,
using Theorem~\ref{th-intersectionord} again, $\CI_r=[f_1,\ldots,f_r]$ is a
reflexive prime difference ideal of dimension $n-r$ and order $\sum_{i=1}^r s_{i}$ over
$\ff\langle \bu_{f_1},\ldots,\bu_{f_r}\rangle$.  Thus, the theorem
is proved. \qedd

\begin{remark}
Notice that Theorem~\ref{inter-thm} is a special case of Theorem~\ref{th-intersectionord},
and the proof of Theorem~\ref{th-intersectionord} gives it an alternative proof.
\end{remark}

\section{Chow form for an irreducible difference variety}

In this section, we will first define the difference Chow form for an irreducible difference variety,
then give its basic properties.

\subsection{Definition of difference Chow form}

Let $V$ be an irreducible difference variety over $\mathcal {F}$ of dimension $d$
and $\mathcal {I}=\mathbb{I}(V)\subset\mathcal {F}\{\mathbb{Y}\}$.
Let
\begin{equation}\label{eq-chowhyp}
\mathbb{P}_i=u_{i0}+u_{i1}y_1+\cdots+u_{in}y_n\,(i=0,\ldots,d)
\end{equation}
be $d+1$ generic difference hyperplanes.
Denote $\bu_i=(u_{i0},u_{i1},\ldots,u_{in})$ for each $i=0,\ldots,d$
and $\bu=\bigcup_{i=0}^n\bu_i\backslash\{u_{i0}\}$.

\begin{lemma}\label{def-Chow}
$[\mathcal {I},\mathbb{P}_0,\ldots,\mathbb{P}_d]\bigcap \mathcal {F}\{\bu_0,\ldots,\bu_d\}$ is a reflexive prime difference ideal of codimension one.
\end{lemma}
\proof
Let $\xi=(\xi_1,\ldots,\xi_n)$ be a generic point of $\mathcal {I}$ which is free from $\mathcal {F}\langle \bu_1,\ldots,\bu_d \rangle$.
Denote $\zeta_i=-\sum_{j=1}^n u_{ij}\xi_j$ and $\zeta=(\zeta_0,u_{01},\ldots,u_{0n},\ldots,\zeta_d,u_{d1},\ldots,u_{dn})$.
We claim that $(\zeta,\xi)$ is a generic point of $\mathcal {J}=[\mathcal {I},\mathbb{P}_0,\ldots,\mathbb{P}_d]\subset\ff\{\bu_0,\ldots,\bu_d,\Y\}$.
It is obvious that $(\zeta,\xi)$ is a zero of $\mathcal {J}$.
Let $g$ be any nonzero difference polynomial in $\mathcal {F}\{\bu_{0},\ldots,\bu_{d},\mathbb{Y}\}$ which vanishes at $(\zeta,\xi)$.
Choose an elimination ranking such that $\bu\prec \Y \prec u_{00} \prec \cdots \prec u_{d0}$.
Then $\mathbb{P}_0,\ldots,\mathbb{P}_d$ constitute an ascending chain.
 Let $r=\prem(g,\mathbb{P}_0,\ldots,\mathbb{P}_d$.
 Then $r\in \mathcal {F}\{\bu,\Y\}$ and $g\equiv r\,\mod\,[\mathbb{P}_0,\ldots,\mathbb{P}_d])$.
 Clearly, $r(\bu,\xi)=0$.
 Since $\bu$ is a set of difference indeterminates over $\ff\langle\xi\rangle$,
  $r\in [\mathcal {I}]\subset \mathcal {F}\{\bu,\Y\}$.
   Hence, $g\in \mathcal {J}$ and it follows that $(\zeta,\xi)$ is a generic zero of $\mathcal {J}$.
Thus $[\mathcal {I},\mathbb{P}_0,\mathbb{P}_1,\ldots,\mathbb{P}_d]\bigcap \mathcal {F}\{\bu_0,\ldots,\bu_d\}$ is a reflexive prime difference ideal with generic zero $\zeta$.

To show that the codimension of $[\mathcal {I},\mathbb{P}_0,\mathbb{P}_1,\ldots,\mathbb{P}_d]\bigcap \mathcal {F}\{\bu_0,\ldots,\bu_d\}$ is $1$,
it suffices to show that $\dtrdeg\,\mathcal {F}\langle \bu\rangle\langle \zeta_0,\ldots,\zeta_d\rangle/\mathcal {F}\langle \bu\rangle=d$.
Since $\dtrdeg\,\mathcal {F}\langle \xi\rangle/\mathcal {F}=d$ and $\zeta_i\in \mathcal {F}\langle \bu,\xi\rangle$,
$\dtrdeg\,\mathcal {F}\langle \bu\rangle\langle \zeta_0,\ldots,\zeta_d\rangle/\mathcal {F}\langle \bu\rangle
\leq  \dtrdeg\,\mathcal {F}\langle \bu\rangle\langle \xi\rangle/\mathcal {F}\langle \bu\rangle=d$.
It is trivial for the case $d=0$.
Now we prove it for the case $d>0$ by showing that  $\zeta_1,\ldots,\zeta_d$ are transformally independent over $\mathcal {F}\langle \bu\rangle$.
Suppose the contrary.  Without loss of generality, assume $\xi_1,\ldots,\xi_d$ is a transformal transcendence basis of $\ff\langle\xi\rangle$ over $\ff$.
Now specialize $u_{ij}$ to $\delta_{ij}(i=1,\ldots,d;j=1,\ldots,n)$,
then by Lemma~\ref{lm-special}, $\xi_{1},\ldots,\xi_{d}$ are transformally dependent over $\mathcal {F}$, which is a contradiction.
Thus $\dtrdeg\mathcal {F}\langle \bu\rangle\langle \zeta_0,\ldots,\zeta_d\rangle/\mathcal {F}\langle \bu\rangle=d$ and the lemma follows.
\qedd

By Lemma~\ref{def-Chow} and Lemma~\ref{le-char-codim1}, there exists a unique irreducible difference polynomial $F(\bu_0,\ldots,\bu_d)$ such  that
$[\mathcal {I},\mathbb{P}_0,\ldots,\mathbb{P}_d]\bigcap \mathcal {F}\{\bu_0,\ldots,\bu_d\}$ is a principal component of $F$\,\footnote{
In differential algebra, it is well known that an irreducible differential polynomial has only one general component.
But in difference case, it is more complicated. In fact, an irreducible difference polynomial $F$ may have more than one principal components
depending on different basic sequences of $F$, which serve as characteristic sets of principal components. For the rigorous definition of {\em basic
sequence}, please refer to \cite{cohn-manifolds}. }.
And from the point view of characteristic sets, if we fix an arbitrary ranking $\mathscr{R}$,
then there exist $F_1,\ldots,F_l$ depending on $\mathscr{R}$ such that
$[\mathcal {I},\mathbb{P}_0,\ldots,\mathbb{P}_d]\bigcap \mathcal {F}\{\bu_0,\ldots,\bu_d\}=\sat(F,F_1,\ldots,F_l)$.

\begin{definition}
The above difference polynomial $F(\bu_0,\ldots,\bu_d)$ is called the {\em difference Chow form} of $V$ or the reflexive prime difference ideal $\mathcal {I}=\mathbb{I}(V)$,
and we call the difference ideal $[\mathcal {I},\mathbb{P}_0,\ldots,\mathbb{P}_d]\bigcap \mathcal {F}\{\bu_0,\ldots, \bu_d\}$ the {\em difference Chow ideal} of $V$.
\end{definition}

The following example shows that the characteristic set of a Chow
ideal could indeed contains more than one elements.
\begin{example}\label{ex-chown=1}
Let $\F=\Q(x)$ and $\sigma(f(x))=f(x+1)$ for each $f\in\Q(x)$. Then $\I=[y_1^2+1,
y_1^{(1)}-y_1]$ is a reflexive prime difference ideal in $\ff\{y_1\}$. The difference Chow form
of $\I$ is $F(\bu_0)=u_{00}^2+u_{01}^2$ and the difference Chow ideal of $\I$ is $\sat(u_{00}^2+u_{01}^2, u_{01}u_{00}^{(1)}-u_{00}u_{01}^{(1)})$.

In general, let $\I = \sat(g(y_1),g_1(y_1),\ldots,g_s(y_1))$ be a
reflexive prime ideal in $\F\{y_1\}$.
Let $F(\bu_0)=M(u_{01})g(-\frac{u_{00}}{u_{01}})$ and $F_i(\bu_0)=M_i(u_{01})g_i(-\frac{u_{00}}{u_{01}})$
where $M(u_{01})$ and $M_i(u_{01})$ are the minimal difference monomials
such that $M(u_{01})g(-\frac{u_{00}}{u_{01}}), M_i(u_{01})g_i(-\frac{u_{00}}{u_{01}})\in\ff\{\bu_0\}$.
Clearly, they are irreducible. The difference Chow ideal of $\I$ is
$\sat\big(F(\bu_0),F_1(\bu_{0}),\ldots,F_s(\bu_{0})\big)$
and the difference Chow form of $\I$ is $F(\bu_0)$.
Indeed, $\mathcal {A}=F,F_1,\ldots,F_s$ constitute an ascending chain in $[\I,u_{00}+u_{01}y_1]\cap\ff\{\bu_0\}$ w.r.t. the elimination ranking $u_{01}\prec u_{00}$.
And if $H$ is a difference polynomial in $[\I,u_{00}+u_{01}y_1]\cap\ff\{\bu_0\}$ which is reduced w.r.t. $\mathcal {A}$, then $H(-u_{01}y_1,u_{01})$ is a difference polynomial in $[\I]\subset\ff\langle u_{01}\rangle\{y_1\}$ reduced w.r.t. $g(y_1),g_1(y_1),\ldots,g_s(y_1)$.
Thus, $H=0$ and $\mathcal {A}$ is a characteristic set of $[\I,u_{00}+u_{01}y_1]\cap\ff\{\bu_0\}$
\end{example}

\begin{example}\label{ex-chowcontrast}
Let $\J=[y_1^2+1, y_1^{(1)}+y_1]$ be a reflexive prime difference ideal in $\Q(x)\{y_1\}$. 
By Example \ref{ex-chown=1},  the difference Chow form of $\J$ is $F(\bu_0)=u_{00}^2+u_{01}^2$ and the difference Chow ideal of $\I$ is $\sat(u_{00}^2+u_{01}^2, u_{01}u_{00}^{(1)}+u_{00}u_{01}^{(1)})$.
Notice that the difference Chow form of $\I=[y_1^2+1, y_1^{(1)}-y_1]$ is equal to that of $\J$.
So different reflexive prime difference ideals may have the same difference Chow form,
which is quite different from the differential case where
the correspondence between differential ideals and differential Chow forms is one-to-one.
Although difference Chow forms can not be used to distinguish different difference ideals,
the correspondence between reflexive prime difference ideals and difference Chow ideals is one-to-one.
So difference Chow ideals play  an important role here.
\end{example}

\begin{example}
Let $\I=[y_1^{(1)}-y_1,y_2^2-y_1,y_2^{(1)}+y_2]\subset\Q\{y_1,y_2\}$.
Then $\I$ is a reflexive prime difference ideal of dimension $0$.
The difference Chow form of $\I$ is $F(\bu_0)=u_{01}u_{02}u_{02}^{(1)}u_{00}^{(1)}+u_{01}u_{00}(u_{02}^{(1)})^2+u_{01}^{2}(u_{00}^{(1)})^2-u_{01}u_{00}u_{01}^{(1)}u_{00}^{(1)}
+u_{02}^2u_{01}^{(1)}u_{00}^{(1)}+u_{00}u_{02}u_{01}^{(1)}u_{02}^{(1)}-u_{01}u_{00}u_{01}^{(1)}u_{00}^{(1)}+u_{00}^{2}(u_{01}^{(1)})^2$
and the difference Chow ideal is $\sat(F(\bu_0),F_1(\bu_0))$ with $F_1(\bu_0)=\left|\begin{array}{ccr}u_{00}&u_{01} & u_{02}\\
u_{00}^{(1)}&u_{01}^{(1)} & -u_{02}^{(1)}\\
u_{00}^{(2)}&u_{01}^{(2)} & u_{02}^{(2)}
\end{array}\right|.$
\end{example}

The following lemma shows that the vanishing of the difference Chow form of $V$ gives a necessary
condition  for a system of difference hyperplanes meeting $V$.
\begin{theorem}
Let $V$ be an irreducible difference variety over $\mathcal {F}$
and $F(\bu_0,\ldots,\bu_d)$ its difference Chow form.
Let $\mathcal {L}_i: a_{i0}+a_{i1}y_1+\cdots+a_{in}y_n=0\,(i=0,\ldots,d)$ be $d+1$ difference hyperplanes
defined over $\mathscr{U}$ and denote $\alpha_i=(a_{i0},\ldots,a_{in})$.
If $V\cap \mathcal {L}_1\cap\cdots\cap \mathcal {L}_d\neq \emptyset$, then the difference Chow ideal of $V$
vanishes at $(\alpha_0,\ldots,\alpha_d)$. In particular, $F(\alpha_0,\ldots,\alpha_d)=0$.
\end{theorem}
\proof Suppose $(\bar{y}_1,\ldots,\bar{y}_n)\in V\cap \mathcal {L}_1\cap\cdots\cap \mathcal {L}_d\neq \emptyset$.
Then the difference ideal $[\mathbb{I}(V),\mathbb{P}_0,\ldots,\mathbb{P}_d]$ vanishes at $(\alpha_0,\ldots,\alpha_d,\bar{y}_1,\ldots,\bar{y}_n)$.
Thus, $[\mathbb{I}(V),\mathbb{P}_0,\ldots,\mathbb{P}_d]\cap \mathcal {F}\{\bu_0,\ldots, \bu_d\}$, in particular $F(\bu_0,\ldots, \bu_d)$, vanishes at $(\alpha_0,\ldots,\alpha_d)$.
\qedd

\begin{remark} \label{remark-computechow}
We remark that the difference characteristic set method proposed in \cite{Gao1}
could be used to compute the difference Chow form of $V$ if we know a set of
finitely many generating  difference polynomials for $V$.
\end{remark}

\subsection{The order of difference Chow form}

In this section, we will show that the order of the difference Chow form is actually equal to that of the corresponding difference variety.
\begin{lemma}\label{chow-ord2}
Let $F(\bu_0,\bu_1,\ldots,\bu_d)$ be the difference Chow form of an irreducible variety $V$ over $\mathcal {F}$.
 Then the following assertions hold.
\begin{itemize}
\item[$1)$] Suppose $F_{\rho\tau}$ is obtained from $F$ by interchanging $\bu_\rho$ and $\bu_\tau$ in $F$. Then $F_{\rho\tau}$ and $F$ differ at most by a sign.

\item[$2)$] $\ord(F,u_{ij})(i=0,\ldots,d;j=0,\ldots,n)$ are the same for all $u_{ij}$ appearing in $F$.
In particular, $u_{i0}$ appears effectively in $F$. And $\ord(F,u_{ij})=-\infty$ if and only if $y_j\in\mathbb{I}(V)$.

\item[$3)$] $\Eord(F,u_{ij})=\ord(F,u_{ij})$, for all the $i,j$.
\end{itemize}
\end{lemma}

\proof
1) Follow the notations in Lemma~\ref{def-Chow}.
Since $\bu$ is a set of difference indeterminates over $\ff\langle\xi\rangle$,
the following difference automorphism $\phi$ of $\mathcal {F}\langle \xi\rangle\langle \bu\rangle$ over $\mathcal {F}\langle \xi\rangle$ can be defined:
$\phi(u_{ij})=u_{ij}^*=\left\{\begin{array}{l}u_{ij}, i\neq \rho,\tau\\ u_{\tau j}, i=\rho\\ u_{\rho j}, i=\tau\end{array}
\right..$ Denote $f(\bu,u_{00},\ldots,u_{d0})=F(\bu_0,\ldots,\bu_d)$, then $f(\bu;\zeta_0,\ldots,\zeta_\rho,\ldots,\zeta_\tau,\ldots,\zeta_d)=0$.
So $\phi\big(f(\bu;\zeta_0,\ldots,\zeta_d)\big)=f(\bu^*;\zeta_0,\ldots,\zeta_\tau,\ldots,\zeta_\rho,\ldots,\zeta_d)=0$.
Let $F_{\rho\tau}(\bu_0,\ldots,\bu_d)=f(\bu^*;u_{00},\ldots,u_{\tau0},\ldots,u_{\rho0},\ldots,u_{d0})$,
then $F_{\rho\tau}(\bu;\zeta_0,\ldots,\zeta_d)=0$.
Thus, $F_{\rho\tau}\in\mathbb{I}(\zeta)=\sat(F,\ldots)$.
Since $\ord(F_{\rho\tau})=\ord(F)$, $\deg(F_{\rho\tau})=\deg(F)$ and $F_{\rho\tau}$ has the same content as $F$,
then $F_{\rho\tau}=\pm F$.

 2) By Lemma~\ref{def-Chow} and 1), we obtain that each $u_{i0}$ appears effectively in $F$ with the same order.
 Suppose $\ord(F,u_{i0})=s$.
 For $j\neq0$, we consider $\ord(F,u_{ij})$.
 If $\ord(F,u_{ij})=l>s$, then we differentiate $f(\bu;\zeta_0,\ldots,\zeta_d)=0$ w.r.t. $u_{ij}^{(l)}$ and we get
 $\frac{\partial f}{\partial u_{ij}^{(l)}}(\bu;\zeta_0,\ldots,\zeta_d)=0$, a contradiction.
 If $\ord(F,u_{ij})=l<s$, differentiate $f(\bu;\zeta_0,\ldots,\zeta_d)=0$ w.r.t. $u_{ij}^{(s)}$,
  then $\frac{\partial f}{\partial u_{i0}^{(s)}}(\bu,\zeta_0,\ldots,\zeta_d)(-\xi_j)=0$.
  Since $\frac{\partial f}{\partial u_{i0}^{(s)}}(\bu,\zeta_0,\ldots,\zeta_d)\neq0$, $\xi_j=0$.
  And $y_j\in \mathbb{I}(V)\Longleftrightarrow \xi_j=0\Longleftrightarrow \zeta_i$ is free of $u_{ij}\Longleftrightarrow F$ is free from $u_{ij}$, thus $\ord(f, u_{ij})=s$ for all $u_{ij}$ appearing in $F$.

 3) Suppose $\lord(F,u_{i0})=t$. Similarly, we can prove that $\lord(F,u_{ij})(i=0,\ldots,d;j=1,\ldots,n)$ are the same for all $u_{ij}$ appearing in F.
 Set $G=\sigma^{(-t)}(F)$, since $\sat(F,\ldots)$ is a reflexive prime difference ideal,
 $G\in \sat(F,\ldots)$. Since $\ord(F)=\ord(G)+t$, $t=0$. Hence $\ord(F,u_{ij})=\Eord(F,u_{ij})$.
\qedd

\begin{definition}
The order of the difference Chow form is defined to be $\ord(F)=\ord(F,u_{i0})$ for any $i \in \{0,\ldots,d\}$.
\end{definition}

The following result shows that the difference characteristic set of
$[\mathcal {I},\P_0,\ldots,\P_d]$ can be easily computed if the
difference characteristic set of the difference Chow ideal is given.
\begin{lemma}\label{chow-ord3}
Let $F(\bu_0,\ldots,\bu_d)$ be the difference Chow form of a
reflexive prime difference ideal $\mathcal {I}$ and
$F,F_1,\ldots,F_l$ a characteristic set of the difference Chow ideal
w.r.t. some ranking $\mathscr{R}$ endowed on $\bigcup_{i=0}^d\bu_i$.
Then
$$\mathcal {A}=\{F,F_1,\ldots,F_l,\frac{\partial F}{\partial u_{00}}y_1-\frac{\partial F}{\partial u_{01}},\ldots,\frac{\partial F}{\partial u_{00}}y_n-\frac{\partial F}{\partial u_{0n}}\}$$
is a characteristic set\footnote{Here $\A$ is a  triangular set but may not be an ascending chain.
Note that the difference remainder of  $\frac{\partial F}{\partial u_{00}}$ is not zero,
so $\A$ can also serve as a characteristic set, which is just similar to the differential case.} of $[\mathcal {I},\P_0,\ldots,\P_d]\subset\mathcal {F}\{\bu_{0},\ldots,\bu_{d},\Y\}$
w.r.t. the elimination ranking $u_{ij}\prec y_1\prec\cdots\prec y_n$ which is consistent with $\mathscr{R}$.
\end{lemma}
\proof Denote $\mathbb{I}_{\zeta,\xi}=[\mathcal {I},\P_0,\ldots,\P_d]\subset\mathcal {F}\{\bu_{0},\ldots,\bu_{d},\Y\}$.
For each $\rho =1,\ldots,n$, differentiate $F(\bu;\zeta_0,\ldots,\zeta_d)=0$ w.r.t. $u_{0\rho}$,
then $\frac{\partial F}{\partial u_{0\rho}}\big|_{(u_{00},\ldots,u_{d0})=(\zeta_0,\ldots,\zeta_d)}-\xi_\rho \frac{\partial F}{\partial u_{00}}\big|_{(u_{00},\ldots,u_{d0})=(\zeta_0,\ldots,\zeta_d)}=0$.
 Hence, $\frac{\partial F}{\partial u_{00}}y_\rho-\frac{\partial F}{\partial u_{0\rho}}\in \mathbb{I}_{\zeta,\xi}(\rho =1,\ldots,n)$.
Let $f$ be any difference polynomial in $\mathbb{I}_{\zeta,\xi}$.
Suppose $g$ is the difference remainder of $f$ w.r.t. $\frac{\partial F}{\partial u_{00}}y_\rho-\frac{\partial F}{\partial u_{0\rho}}\,(\rho=1,\ldots,n)$ w.r.t. the elimination ranking $u_{ij}\prec y_1\prec\cdots\prec y_n$,
then $g\in\mathbb{I}_{\zeta,\xi}\cap\ff\{\bu_0,\ldots,\bu_d\}$.
Thus, $\prem(f,\mathcal{A})=\prem(g,[F,F_1,\ldots,F_l])=0$.
Therefore $\mathcal {A}$ is a characteristic set of $\mathbb{I}_{\zeta,\xi}$
w.r.t. the elimination ranking $u_{ij}\prec y_1\prec\cdots\prec y_n$ which is consistent with $\mathscr{R}$.
\qedd

The following result shows that the generic point $(\xi_1,\ldots,\xi_n)$
of $V$ can be recovered from its difference Chow ideal.

\begin{cor}
Let $F(\bu_{0},\ldots,\bu_{d})$ be the difference Chow form of $V$.
Suppose $\zeta$ is a generic point of the difference Chow ideal of $V$ and denote
\[\eta_{\rho}=\overline{\frac{\partial F}{\partial u_{0\rho}}}\bigg/\overline{\frac{\partial F}{\partial u_{00}}}\,\,\,(\rho=1,\ldots,n)\]
where $\overline{\frac{\partial F}{\partial u_{0\rho}}}=\frac{\partial f}{\partial u_{0\rho}}\big|_{(\bu_{0},\ldots,\bu_{d})=\zeta}$.
Then $(\eta_1,\ldots,\eta_n)$ is a generic point of $V$.
\end{cor}
\proof It follows directly from
Lemma~\ref{chow-ord3}. \qedd

\begin{cor}
Let $V$ be an irreducible difference variety of dimension $d$ over
$\ff$. Denote $\bu=\cup_{i=0}^d\bu_i\backslash\{u_{i0}\}$. Then over
$\ff\langle\bu\rangle$, $V$ is birationally equivalent to an
irreducible difference variety of codimension $1$.
\end{cor}
\proof Suppose $F(\bu_0,\ldots,\bu_d)$ is the difference Chow form
of $V$ and $\text{\rm CI}\subset\ff\{\bu;u_{00},\ldots,u_{d0}\}$ is
the difference Chow ideal of $V$. Let $\text{\rm
CI}_{\bu}=[\text{\rm
CI}]\subset\ff\langle\bu\rangle\{u_{00},\ldots,u_{d0}\}$ and
$W=\V(\text{\rm CI}_{\bu})\subset\mathscr{U}^n$. By
Lemma~\ref{ord1}, $W$ is an irreducible difference variety of
codimension $1$. Then over $\ff\langle\bu\rangle$, $V$ is
birationally equivalent to $W$ with the following maps:
\[\begin{array}{lccc}
\phi:& V & \longrightarrow & W\\
\quad&(a_1,\ldots,a_n)&\quad&
(-\sum\limits_{k=0}^nu_{0k}a_k,\ldots,-\sum\limits_{k=0}^nu_{dk}a_k)\\
\end{array}\]
and \[\begin{array}{lccc}
\psi:& W & \longrightarrow & V\\
\quad&(b_{00},\ldots,b_{d0})&\quad& (\overline{\frac{\partial
F}{\partial u_{01}}}\big/\overline{\frac{\partial F}{\partial u_{00}}},
\ldots,\overline{\frac{\partial F}{\partial u_{0n}}}\big/\overline{\frac{\partial F}{\partial u_{00}}}),\\
\end{array}\]
where $\overline{\frac{\partial F}{\partial u_{0k}}}=\frac{\partial
F}{\partial u_{0k}}(\bu;b_{00},\ldots,b_{d0})$. \qedd

The following result gives our first main property for difference Chow form.
\begin{theorem}\label{chow-ord-thm}
Let $\mathcal {I}$ be a reflexive prime difference ideal of dimension $d$
with  difference Chow form $F(\bu_{0},\ldots,\bu_{d})$.
Then $\ord(F)=\ord(\mathcal {I})$.
\end{theorem}

\proof
Let $\mathcal {I}_d=[\mathcal {I},\mathbb{P}_1,\ldots,\mathbb{P}_d]\subset\mathcal {F}\langle\bu_1,\ldots,\bu_d\rangle\{\mathbb{Y}\}$.
 By Theorem~\ref{ord-thm} $\mathcal {I}_d$ is a reflexive prime difference ideal with $\dim(\mathcal {I}_d)=0$ and $\ord(\mathcal {I}_d)=\ord(\mathcal {I})$.

let $\mathcal {J}=[\mathcal {I},\mathbb{P}_0,\ldots,\mathbb{P}_d]=[\mathcal {I}_d,\mathbb{P}_0]
\subset\mathcal {F}\langle \bu_1,\ldots,\bu_d;u_{01},\ldots,u_{0n}\rangle\{u_{00},\mathbb{Y}\}$.
Choose a ranking $\mathscr{R}$ such that $u_{00}$ is the leading variable of $F$.
By Lemma~\ref{chow-ord3}, $\mathcal {A}$ is a characteristic set of $\mathbb{I}_{\xi,\zeta}$.
Since $\{\bu_1,\ldots,\bu_d,u_{00},\ldots,u_{0n}\}$ is a parametric set of $\mathbb{I}_{\xi,\zeta}$,
 by Lemma~\ref{ord1}, $\mathcal {A}$ is also a characteristic set of $\mathcal {J}$ w.r.t some ranking.
Since $\dim(\mathcal {J})=0$, $\ord(\mathcal {J})=\ord(\mathcal {A})=\ord(F)$.

 Let $\eta=(\eta_1,\ldots,\eta_n)$ be a generic zero of $\mathcal {I}_d$.
 Set $\theta=-\sum_{j=1}^nu_{0j}\eta_j$,
 then $(\theta,\eta_1,\ldots,\eta_n)$ is a generic zero of $\mathcal {J}$.
 Since $\dim(\mathcal {J})=0$,  for sufficiently large integer $t$,
 \begin{eqnarray}
\ord(\mathcal {J})
&=&\varphi_{\mathcal {J}}(t) \nonumber \\
&=&\trdeg\,\mathcal {F}\langle \bu_1,\ldots,\bu_d,u_{01},\ldots,u_{0n}\rangle(\theta^{[t]},\eta_1^{[t]},\ldots,\eta_n^{[t]})/\mathcal {F}\langle \bu_1,\ldots,\bu_d,u_{01},\ldots,u_{0n}\rangle \nonumber \\
&=&\trdeg\,\mathcal {F}\langle \bu_1,\ldots,\bu_d,u_{01},\ldots,u_{0n}\rangle(\eta_1^{[t]},\ldots,\eta_n^{[t]})/\mathcal {F}\langle \bu_1,\ldots,\bu_d,u_{01},\ldots,u_{0n}\rangle \nonumber \\
&=&\trdeg\,\mathcal {F}\langle \bu_1,\ldots,\bu_d\rangle(\eta_1^{[t]},\ldots,\eta_n^{[t]})/\mathcal {F}\langle \bu_1,\ldots,\bu_d\rangle \nonumber \\
&=&\varphi_{\mathcal {I}_d}(t)=\ord(\mathcal {I}_d) \nonumber
\end{eqnarray} Hence $\ord(F)=\ord(\mathcal {J})=\ord(\mathcal {I}_d)=\ord(\mathcal {I})$.
\qedd

\subsection{Homogeneity of the difference Chow form}
 In this section, we will show that the difference Chow form is transformally homogenous.

\begin{definition}
A difference polynomial $p \in \mathcal {F}\{y_0,\ldots,y_n\}$ is said to be  transformally homogenous
if for a new difference indeterminate $\lambda$,
$p(\lambda y_0,\ldots,\lambda y_n)=M(\lambda)p(y_0,\ldots,y_n)$,
 where $M(\lambda)$ is a difference monomial of $\lambda$.
\end{definition}

The difference analog of Euler's theorem related to homogeneous polynomials is valid.
\begin{lemma}\label{chow-homo1}\cite{li-sddr}
A difference polynomial $p \in \mathcal {F}\{y_0,\ldots,y_n\}$ is
transformally homogeneous if and only if for each $k\geq0$, there
exists $r_k\in\mathbb{N}_0$ such that
$$\sum_{j=0}^ny_j^{(k)}\frac{\partial p}{\partial y_j^{(k)}}=r_kp.$$
\end{lemma}

\begin{theorem}\label{chow-homo-thm}
Let $F(\bu_0,\bu_1,\ldots,\bu_d)$ be the difference Chow form of a difference irreducible variety $V$ of dimension $d$ and order $h$. Then

$1)$ $\sum_{j=0}^n u_{\tau j}^{(k)}\frac{\partial F}{\partial u_{\sigma j}^{(k)}}=\left\{
\begin{array}{ll}
0& \text{if}\,\sigma=\tau\\
r_kF& \text{if}\,\sigma\neq \tau
\end{array}
\right.$
for $k=0,1,\ldots,h$, where $r_k\in\mathbb{N}_0$.

$2)$ $F(\bu_0,\ldots\bu_d)$ is transformally homogeneous  in each $\bu_i$.
\end{theorem}

\proof
Differentiate  $F(\bu;\zeta_0,\zeta_1,\ldots,\zeta_d)=0$ on both sides w.r.t. $u_{\sigma j}^{(k)}(k=0,\ldots,h)$, then
\begin{equation}
\overline{\frac{\partial F}{\partial u_{\sigma j}^{(k)}}}+\overline{\frac{\partial F}{\partial u_{\sigma 0}^{(k)}}}(-\xi_j^{(k)})=0, \nonumber
\end{equation}
Where $\overline{\frac{\partial F}{\partial u_{\sigma j}^{(k)}}}=\frac{\partial F}{\partial u_{\sigma j}^{(k)}}\big|_{(u_{00},\ldots,u_{d0})=(\zeta_0,\ldots,\zeta_d)}$.
Multiply the above equation by $u_{\tau j}^{(k)}$ and for $j$ from $1$ to $n$,
add them together, then we get
\begin{equation}
\sum_{j=1}^nu_{\tau j}^{(k)}\overline{\frac{\partial F}{\partial u_{\sigma j}^{(k)}}}+\zeta_\tau^{(k)}\overline{\frac{\partial F}{\partial u_{\sigma 0}^{(k)}}}=0. \nonumber
\end{equation}

Hence, the difference polynomial $G_{\sigma\tau}=\sum_{j=0}^nu_{\tau j}^{(k)}\frac{\partial F}{\partial u_{\sigma j}^{(k)}}\in\mathbb{I}_{\zeta,\xi}$.
Since $\ord(G_{\sigma\tau})\leq h$,
$F$ divides $G_{\sigma\tau}$.
If $\tau\neq \sigma$, $\deg(G_{\sigma\tau},\bu_\sigma^{(k)})<\deg(F,\bu_\sigma^{(k)})$, thus $G_{\sigma\tau}=0$.
In the case $\tau=\sigma$, there exists $r_k\in\mathbb{N}_0$ such that
$\sum_{j=0}^nu_{\sigma j}^{(k)}\frac{\partial F}{\partial u_{\sigma j}^{(k)}}=r_kF$ for $k=0,1,\ldots,h$.
And by  Lemma~\ref{chow-homo1}, $F(\bu_0,\bu_1,\ldots\bu_d)$ is transformally homogeneous  in each $\bu_i$.
\qedd

\begin{definition}\label{def-ddeg}
Let $V$ be an irreducible difference variety of dimension $d$ and order $h$.
Let $F(\bu_{0},\ldots,\bu_{d})$ be the difference Chow form of $V$.
The {\em difference degree} of $V$ is defined as the homogenous degree $r=\sum_{k=0}^hr_k$ of its difference Chow form in each
$\bu_{i}$ $(i=0,\ldots,d)$.
\end{definition}

The following result shows that the difference degree of a variety $V$
is an invariant of $V$ under invertible linear transformations.

\begin{lemma}\label{lm-cfl}
Let $A=(a_{ij})$ be an $n\times n$ invertible matrix with
$a_{ij}\in\mathcal {F}$ and $F(\bu_0,\bu_1,\ldots,\bu_d)$ the Chow
form of an irreducible difference variety $V$ of dimension $d$. Then
the difference Chow form of the image variety of $V$ under the linear
transformation $\Y=A\X$ is $F^{A}(\bv_{0},\ldots,\bv_{d})=F(\bv_{0}B,\ldots,\bv_{d}B)$, where
$B=\left(\begin{array}{ll} 1 &  0_{1\times n}  \\0_{n\times 1}  &A \\
\end{array}\right)$ and $\bu_{i}$ and $\bv_{i}$ are
regarded as row vectors.
\end{lemma}
\proof Let $\xi=(\xi_{1},\ldots,\xi_{n})$ be a generic point of $V$.
Under the linear transformation $\Y=A\X$, $V$ is mapped to an irreducible difference variety $V^A$ whose generic point is $\eta=(\eta_{1},\ldots,\eta_{n})$ with $\eta_i=\sum_{j=1}^n
a_{ij}\xi_j$. Denote $F(\bu_0,\ldots,\bu_d)=f(u_{ij};u_{00},\ldots,\bu_d)$.
Note that $F^{A}(\bv_{0},\ldots,\bv_{d})=f(\sum_{k=1}^nv_{ik}a_{kj};v_{00},\ldots,v_{d0})$
 and $f(\sum_{k=1}^n v_{ik}a_{kj};$ $-\sum_{k=1}^n v_{0k}\eta_{k},\ldots,$ $-\sum_{k=1}^n v_{dk}\eta_{k})
 $ $=f(\sum_{k=1}^n v_{ik}a_{kj};-\sum_{j=1}^n(\sum_{k=1}^n
v_{0k}a_{kj})\xi_{j},\ldots,-\sum_{j=1}^n(\sum_{k=1}^n v_{dk}a_{kj})\xi_{j})=0.$ Since $V^A$ is of the same dimension
and order as $V$  and  $F^{A}$ is irreducible, by the definition
of difference Chow form, the proof is completed.\qedd

\begin{definition} \label{denomination}
Let  $p$ be a difference polynomial in $\mathcal {F}\{y\}$.
Suppose $\ord(p,y)=t$ and $m_i=\deg(p,y^{(i)})\,(i=0,\ldots,t).$
Then $\prod_{i=0}^t(y^{(i)})^{d_i}$ is called the {\em difference denomination} of $p$,
denoted by $\den^\sigma(p)$.
\end{definition}

\begin{example}\label{ex-chowdegree}
 Consider the case $n=1$.
Suppose $\CI=\sat(g(y),g_1(y),\ldots,g_s(y))$ be a reflexive prime difference ideal in $\ff\{\Y\}$.
Let $\den^\sigma(g)=M(y).$
Clearly, $M(u_{01})$ is the minimal difference monomial such that  $M(u_{01})g(-\frac{u_{00}}{u_{01}})\in\ff\{\bu_0\}$ where $\bu_0=(u_{00},u_{01})$.
By Example~\ref{ex-chown=1}, $F(\bu_0)=M(u_{01})g(-\frac{u_{00}}{u_{01}})$ is the difference Chow form of $\CI$.
Thus, the difference degree of $\CI$ is equal to the degree of the difference denomination of $g$, i.e. $\sum_{i=0}^{\ord(g)}\deg(g,y^{(i)})$.
Recall that in the differential case, the differential degree of $\sat(g)$ is also equal to the denomination of $g.$
But the denomination of a differential polynomial is much more complicated to compute than the difference case.
There, we showed that the weighted degree of $g$ is a sharp bound for the differential  degree of $\sat(g)$.
\end{example}

\section{Conclusion}

In this paper, firstly, it is  shown that both the dimension and the order
of a reflexive prime difference ideal can be read off from its
characteristic sets under any fixed ranking.
Then we give a generic
intersection theorem for difference varieties. Precisely, the
intersection of an irreducible difference variety of dimension $d >
0$ and order $h$ with a generic difference hypersurface of order $s$ is shown to be
an irreducible difference variety of dimension $d-1$ and order $h+s$.
Based on the intersection theory, the difference Chow form for an
irreducible difference variety is defined and  its basic properties
are given.

Below, we propose several problems for further study.

In the differential case, much more properties are proved for the differential Chow form \cite{gao-dcf},
which are not yet able to be generalized to the difference case due to the distinct structures of the differential and
difference operators. It is interesting to enrich the properties of difference Chow form, especially to establish
a theory of difference Chow variety.

In Remark \ref{remark-computechow}, we mentioned that the difference Chow form can be computed with the difference characteristic set method.
But it is difficult to analyze the computing complexity if we just work with the usual characteristic set method.
In the algebraic case, Jeronimo et al gave a bounded probabilistic algorithm to compute
the Chow form, whose complexity is polynomial in the size and the geometric degree of the
input equation system \cite{Complexitychowform}.
It is important to apply the principles behind  such algorithms to propose an efficient algorithm to compute the difference Chow form.

In Theorem~\ref{ord-thm},  we proved that both the dimension and the relative order of a reflexive prime difference ideal can be
reflected from its characteristic set under an arbitrary ranking.
We conjecture that the relative effective order can also be read off.

\end{document}